\documentclass[11pt]{amsart}
\usepackage{amsmath,amsfonts,amsthm,amssymb,amsthm,graphicx,mathtools,amscd,enumerate,paralist,bm,hyperref}

\usepackage[mathscr]{eucal}
\usepackage{graphics,color,epsfig}
\usepackage[nocompress]{cite}

\DeclareSymbolFont{iwonaletters}{OML}{iwona}{m}{it}
\DeclareMathSymbol{\bdel}{\mathalpha}{iwonaletters}{"E}
\DeclareMathOperator*{\argmax}{arg\,max}

\theoremstyle{plain}
\newtheorem{theorem}{Theorem}[section]
\newtheorem{lemma}[theorem]{Lemma}

\newtheorem{definition}[theorem]{Definition}
\newtheorem{example}[theorem]{Example}
\newtheorem{proposition}[theorem]{Proposition}

\newtheorem{assumption}[theorem]{Assumption}

\newtheorem{problem}[theorem]{Problem}
\newtheorem{remark}[theorem]{Remark}

\newcommand{\la}{\lambda}
\newcommand{\eps}{\varepsilon}
\newcommand{\ph}{\varphi}

\newcommand{\gam}{\gamma}
\newcommand{\kap}{\kappa}
\newcommand{\sig}{\sigma}
\newcommand{\del}{\delta}
\newcommand{\om}{\omega}

\newcommand{\Del}{\mathnormal{\Delta}}

\newcommand{\N}{{\mathbb N}}
\newcommand{\Q}{{\mathbb Q}}
\newcommand{\R}{{\mathbb R}}
\newcommand{\U}{{\mathbb U}}

\newcommand{\Z}{{\mathbb Z}}

\newcommand{\E}{{\mathbb E}}

\newcommand{\PP}{{\mathbb P}}

\newcommand{\calB}{{\mathcal B}}

\newcommand{\calL}{{\mathcal L}}
\newcommand{\calM}{{\mathcal M}}

\newcommand{\calS}{{\mathcal S}}

\newcommand{\bu}{\mathbf{u}}
\newcommand{\bv}{\mathbf{v}}
\newcommand{\bw}{\mathbf{w}}

\newcommand{\bB}{\mathbf{B}}

\newcommand{\bLa}{\mathbf{\Lambda}}

\newcommand{\frp}{\mathfrak{p}}
\newcommand{\frs}{\mathfrak{s}}

\renewcommand{\proof}{\noindent{\bf Proof.\ }}

\newcommand{\lan}{\langle}
\newcommand{\ran}{\rangle}

\newcommand{\w}{\wedge}

\newcommand{\pl}{\partial}

\newcommand{\To}{\Rightarrow}

\newcommand{\grad}{\nabla}
\newcommand{\iy}{\infty}

\newcommand{\loc}{{\rm loc}}

\newcommand{\noi}{\noindent}

\newcommand{\Ll}{\mathbb{L}}

\begin{document}

\title[BRW and free obstacle]{
Trimmed branching random walk\\ and a free obstacle problem
}

\author{Rami Atar}
\address{Viterbi Faculty of Electrical \& Computer Engineering
\\
Technion
} 
\email{rami@technion.ac.il}

\author{Leonid Mytnik}
\address{Faculty of Data \& Decision Sciences
\\
Technion
} 
\email{leonidm@technion.ac.il}

\author{Gershon Wolansky}
\address{Faculty of Mathematics
\\
Technion
}
\email{gershonw@technion.ac.il}


\date{\today}

\begin{abstract}
Consider $N$ particles performing random walks on the $\eps$-grid $(\eps\Z)^d$, $\eps>0$ with branching and density-dependent selection: When one of the particles branches, a particle is removed from the most populated site. The walks are assumed to be asymptotic, as $\eps\to0$, to diffusion processes of the form
\[
dX_i(t)=b(X_i(t))dt+\sqrt{2}dW_i(t),
\]
for $b$ a given vector field. Denoting $\calL^*=\Del-\grad\cdot(b\,\cdot)$, the hydrodynamic limit, as $N\to\iy$ followed by $\eps\to0$, is characterized in terms of a parabolic free obstacle problem
\[
\pl_tu=\calL^*u+u-\beta
\]
where $\beta$ is a measure on $\R^d\times[0,\iy)$ supported on $\{(x,t):u(x,t)=|u(\cdot,t)|_\iy\}$. Here, the unknowns are $u$, the mass density, and $\beta$, the removal measure, for which $t\mapsto\beta(\R^d\times[0,t])$ is prescribed. This is analogous to the well-understood relation between particle systems with spatial selection and free boundary problems, but the techniques require quite different ideas. The key ingredients of the proof include PDE uniqueness for continuous densities and a uniform-in-$\eps$ estimate on modulus of continuity of prelimit densities. The work gives rise to open problems such as ``flat top'' versus ``sharp top'' solutions, which are discussed based on concrete examples.
\end{abstract}

\maketitle

\section{Introduction}

\subsection{Setting and main result}
We consider $N$ continuous time random walks on the
$\eps$-grid $\calS_\eps:=(\eps\Z)^d$, $0<\eps<1$, with branching and selection.
Branching occurs according to a
rate-$1$ Poisson process for each particle. Upon branching, one particle
is removed from the most populated site. As $\eps\to0$, the walks are asymptotic to mutually independent diffusion processes given by
\[
dX_i(t)=b(X_i(t))dt+\sqrt{2}dW_i(t),
\]
where $b$ is a given vector field. We are interested in characterizing the hydrodynamic limit, as $N\to\iy$ followed by $\eps\to0$, in terms of a free obstacle problem. In a scenario where each site can contain no more than $hN$ particles, $h>0$ a given constant, and particles are removed when this number is exceeded, the limiting particle density $u$ will satisfy a parabolic PDE with an obstacle given by $u(x,t)\le h$. In our setting, trimming is applied in a way that the number of particles is preserved, leading to a global condition $\int u(x,t)dx=1$ and an obstacle of the form $u(x,t)\le h(t)$, where $h$ is free (i.e., not prescribed).

To give a precise description of the setting, let $k_i$, $1\le i\le 2d$ be defined as
$k_i=e_i$, $k_{i+d}=-e_i$, $1\le i\le d$,
where $e_i$ are the coordinate vectors in $\R^d$.
The generator of the single-particle Markov process is given by
\begin{equation}\label{22a}
\calL_\eps f(x)=\sum_{i\le2d}r_{\eps,i}(x)(f(x+\eps k_i)-f(x)), \qquad x\in\calS_\eps,
\end{equation}
where the jump rates have the form
\[
r_{\eps,i}(x)=\eps^{-2}+\eps^{-1}q_{\eps,i}(x),
\qquad x\in\calS_\eps,\ 1\le i\le 2d.
\]
Denote
\begin{equation}\label{22b}
b_\eps(x)=\sum_{i\le2d}q_{\eps,i}(x)k_i,
\qquad
x\in\calS_\eps.
\end{equation}
Throughout, $|\cdot|$ denotes Euclidean norm in $\R^d$.
\begin{assumption}[Drift]
\label{assn1}
\ \\ (i)
There exists a vector field $b\in C^1_b(\R^d,\R^d)$ with $b$ and $\nabla b$ globally Lipschitz, such that with $b_\eps$ defined via \eqref{22b},
\begin{equation}\label{c0}
\lim_{\eps\to0}
\sup_{x\in\calS_\eps}|b_\eps(x) - b(x)|=0.
\end{equation}
(ii)
There exist $C_1\in(0,\iy)$ and $\eps_0\in(0,1)$
such that for all $\eps\in(0,\eps_0)$, $x,y\in\calS_\eps$
and $i\in[2d]$, one has $q_{\eps,i}(x)\le C_1$ and
\begin{equation}\label{c01}
|q_{\eps,i}(x)-q_{\eps,i}(y))|
\le C_1|x-y|, \qquad
\Big|\frac{q_{\eps,i}(x+\eps k_i)-q_{\eps,i}(x)}{\eps}
-\frac{q_{\eps,i}(y+\eps k_i)-q_{\eps,i}(y)}{\eps}\Big|
\le C_1|x-y|.
\end{equation}
\end{assumption}

As we note in Appendix \ref{app1}, given $b\in C^1_b(\R^d,\R^d)$ with $b$ and $\nabla b$ globally Lipschitz, there always exists $q_\eps$ such that $(q_\eps,b_\eps,b)$ satisfy Assumption \ref{assn1}.

Let the collection of living particles at time $t$ be denoted
by $X_i(t)=X^N_i(t)$, $i\in\{1,\ldots,N\}$. Let
\[
\xi^N_t(dx)=\frac{1}{N}\sum_i\del_{X_i(t)}(dx).
\]
To emphasize the dependence on $\eps$
we will sometimes write the above as $\xi^{\eps,N}_t$.

For $p\in[1,\iy]$, let $|f|_p$ denote the $\Ll_p$ norm of a function $f:\R^d\to\R$ or $f:\calS_\eps\to\R$. Let the space of finite (respectively, probability) Borel measures on $\R^d$ be denoted by $\calM(\R^d)$ (respectively, $\calM_1(\R^d)$), and equip it with the topology of weak convergence. Following is an assumption on the initial empirical measures, $\xi^N_0$.
\begin{assumption}[Initial condition]
\label{assn2}
There exists a uniformly continuous $u_0\in C(\R^d,\R_+)$, and for every $0<\eps<1$ there exists $\bu_0^{(\eps)}:\calS_\eps\to\R_+$
such that the follows holds.
\\
i. $|u_0|_1=1$ and $|\bu_0^{(\eps)}|_1=1$.
\\
ii.
Denoting
\[
\xi_0^{(\eps)}=\sum_{x\in\calS_\eps}\bu^{(\eps)}_0(x)\del_x,
\]
one has $\xi^{\eps,N}_0\to\xi^{(\eps)}_0$ in probability in
$\calM_1(\R^d)$ as $N\to\iy$.
\\
iii. One has
$\sup_{x\in\calS_\eps}|\eps^{-d}\bu_0^{(\eps)}(x)-u_0(x)|\to0$
as $\eps\to0$.
\end{assumption}
\begin{remark}
Consider $Y_i^N$, $i=1,\ldots,N$ i.i.d.\ with a uniformly continuous density $u_0$
and let $X^{\eps,N}_i(0)$, $i=1,\ldots,N$ be a discretization of $Y_i$ such as
$\eps\lfloor\eps^{-1}Y_i \rfloor$, with $\lfloor\cdot\rfloor$ applied componentwise.
Then Assumption \ref{assn2} holds with $u_0^{(\eps)}(x)=\int_{[x_1,x_1+\eps]\times\cdots\times[x_d,x_d+\eps]} u_0(y)dy$ for $x\in\calS_\eps$.
\end{remark}

We are interested in a macroscopic description of the particle system in terms of a PDE.
For $f\in C^2(\R^d,\R)$, denote
\[
\calL f=\Del f+b\cdot\nabla f,
\qquad
\calL^*f=\Del f-\nabla\cdot(b f).
\]
The PDE is concerned with a pair $(u,\beta)$, where
$u(\cdot,t)$ is the macroscopic density at time $t$,
and $\beta$ describes the distribution of mass removal. 
Let $\calM_\loc(\R^d\times\R_+)$ denote the space of
Borel measures on $\R^d\times\R_+$ that are finite
on $\R^d\times[0,T]$ for every $T$, equipped with
the topology of weak convergence on $\R^d\times[0,T]$
for every $T$. Let
\[
\calM^{(1)}(\R^d\times\R_+)=\{
\mu\in\calM_\loc(\R^d\times\R_+):\mu(\R^d\times[0,T])=T,
\,T>0\}.
\]
Any $\mu\in\calM^{(1)}(\R^d\times\R_+)$ can be
disintegrated as $\mu(dx,dt)=\mu_t(dx)dt$
where $\mu_t$ is a probability measure for a.e.\ $t$, and the notation $\mu_t$ will be used for this disintegration, throughout the paper.
Let $\calM^{(1)}(\calS_\eps\times\R_+)$ be defined analogously.

For $u\in C_b(\R^d\times\R_+,\R_+)$ let 
\begin{equation}\label{betau} \bB^u:= \left\{\beta\in \calM^{(1)}(\R^d\times\R_+)  \ : \ \beta_t(\{x:u(x,t)<|u(\cdot,t)|_\iy\})=0
\text{ for a.e.\ $t$} \right\}.
\end{equation}
For $(u,\beta)\in C_b(\R^d\times\R_+,\R_+)\times
\calM^{(1)}(\R^d\times\R_+)$, the PDE reads
\begin{equation}\label{c1}
\pl_t u= \calL^* u+u-\beta,
\qquad u(\cdot,0)=u_0, \qquad \beta\in\bB^u \ . 
\end{equation}

\begin{definition}[Solution to \eqref{c1}]
\label{def1}
A solution to \eqref{c1} is a pair
$(u,\beta)\in C_b(\R^d\times\R_+,\R_+)\times\calM^{(1)}(\R^d\times\R_+)$ satisfying $\beta\in\bB^u$  and for any $t\geq 0$ and any  $\ph\in C^\iy_c(\R^d)$,
\[
\lan \ph,u(\cdot,t)\ran=\lan\ph,u_0\ran
+\int_0^t\lan \calL\ph+\ph,u(\cdot,s)\ran ds
-\int_0^t\lan \ph,\beta_s\ran ds \ . 
\]
\end{definition} 

We now construct the processes $\xi^N_t$ via Poisson random measures (PRM).
Fix $\eps$. For the purpose of tie breaking, fix a total order $\le_{\rm to}$
on the grid $\calS_\eps$.
For a finite nonzero measure $\xi$ on $\calS_\eps$,
let $\argmax\xi=\{x:\xi(x)=\max_y\xi(y)\}$
 and $\argmax^*\xi$ be the unique $x\in\argmax\xi$ such that
 $x\le_{\rm to}y$ for all $y\in\argmax\xi$. Define similarly $\argmax^*u$ for a  function $u:\calS_\eps\to\R_+$ with $0<|u|_1<\iy$.

Let $\Pi$ and $\Gamma$ be PRM on $\calS_\eps\times[2d]\times\R_+\times\R_+$ and, respectively, $\calS_\eps\times\R_+\times\R_+$,
with intensity measures
$(\sum_{y\in\calS_\eps}\del_y(dx))(\sum_{j\le 2d}\del_j(di))d\theta\,dt$ and respectively,
$(\sum_{y\in\calS_\eps}\del_y(dx))d\theta\,dt$.
Assume that the initial condition $\xi^N_0$, $\Pi$ and $\Gamma$
are mutually independent. Then $\xi^N$ is given by
\begin{align}\notag
\xi^N_t=\xi^N_0
&+\frac{1}{N}\int_{\calS_\eps\times[2d]\times\R_+\times[0,t]}1_{[0,Nr_i(x)\xi^N_{s-}(x)]}(\theta)
(\del_{x+\eps k_i}-\del_x)
\Pi(dx,di,d\theta,ds)
\\ \label{20a}
&+\frac{1}{N}\int_{\calS_\eps\times\R_+\times[0,t]} 1_{[0,N\xi^N_{s-}(x)]}(\theta)
(\del_x-\del_{\argmax^*\xi^N_{s-}})
\Gamma(dx,d\theta,ds).
\end{align}
The process $\xi^N_t$ has sample paths in $D(\R_+,\calM(\R^d))$. We equip $D$ with the induced Skorohod topology. Define
\[
\beta^N(dx,dt)=\del_{\argmax^*\xi^N_t}(dx)dt.
\]
Let $\U(\R^d\times\R_+)$ denote the set of functions $u\in C(\R^d\times\R_+,\R_+)$ satisfying, for every $T$, $\sup_{t\in[0,T]}|u(\cdot,t)|_\iy<\iy$, and the condition
\[
\text{$x\mapsto u(x,t)$ is continuous in $x$ uniformly in $(x,t)\in\R^d\times[0,T]$.}
\]
Our main result is as follows.

\begin{theorem}\label{th1}
Let Assumptions \ref{assn1} and \ref{assn2} hold.
\\
(a) Consider equation \eqref{c1} with initial condition $u_0$ as in Assumption \ref{assn2}. Then, within the class  $\U(\R^d\times\R_+)\times\calM^{(1)}(\R^d\times\R_+)$, there exists a unique solution $(u,\beta)$ to \eqref{c1}.
\\
(b) Let
\[
\xi_0(dx)=u_0(x)dx, \qquad \xi_t(dx)=u(x,t)dx,\ t>0.
\]
Then there exists a sequence $\eps_N\downarrow0$ such that,
with $(\xi^{(N)},\beta^{(N)})=(\xi^{\eps_N,N},\beta^{\eps_N,N})$, one has
$(\xi^{(N)},\beta^{(N)})\to(\xi,\beta)$
in $D(\R_+,\calM(\R^d))\times\calM^{(1)}(\R^d\times\R_+)$
in probability. 
\end{theorem}

\subsection{Related work}

Particle systems with selection were proposed in \cite{bru06, bru07} as models for natural selection in population dynamics, where the position of a particle on the real line represents the degree of fitness of an individual to its environment. A closely related model,
referred to as the $N$-particle branching Brownian motion ($N$-BBM), was introduced in \cite{mai16}. Here, $N$ Brownian particles on the line are subject to branching and selection, where upon each branching event, the leftmost particle is removed.
At the hydrodynamic limit, it was shown to give rise to a free boundary problem (FBP) \cite{de-masi-nbbm, ber19}.
Various other particle systems with selection, with and without branching, were characterized by FBP at the hydrodynamic limit, including Brownian systems in dimension one \cite{de-masi-book, ata25} and higher \cite{ber21, ber22bee}, and non-local branching models \cite{dur11, de-masi-nbbm2, ata-dr}.
Apart from models involving selection, work on particle systems that are macroscopically described by FBP include \cite{de-masi-excl} which studies a variant of the simple exclusion process, and \cite{dem19, ata-bud-25} which study the Atlas model.

A population model with selection was studied in \cite{fou-mel}, where individuals die at rate that depends on the local population density, establishing macroscopic approximations. Our model, like \cite{fou-mel}, captures density dependent selection. It also exhibits clear similarities to the $N$-BBM, where removals occur at the extremes: specifically, the leftmost particle in the $N$-BBM and the most populated site in our model.

We now draw an analogy between the macroscopic model that arises from the $N$-BBM \cite{de-masi-nbbm, ber19}, and the one obtained in our work. In the former case, this is a FBP that takes the following form. Given an initial density $u_0$, $|u_0|_1=1$, find $(u,\sig)$ solving
\begin{equation}\label{n1}
\begin{cases}
\pl_t=\pl_{xx} u+u & x>\sig(t),
\\
u=0 & x\le \sig(t),
\\
|u(\cdot,t)|_1=1 & t>0,
\\
u(\cdot,0)=u_0.
\end{cases}
\end{equation}
Here, $u$ represents the mass density and $\sig$ the left edge of its support.
Consider now equation \eqref{c1} as a free obstacle problem. Given $u_0$, $|u_0|_1=1$, find $(u,h)$ such that
\begin{equation}\label{29a}
\begin{cases}
0\le u(x,t)\le h(t)& x\in\R^d,t>0,
\\
\pl_tu(x,t)=\calL^*u(x,t)+u(x,t)& \text{\rm if } u(x,t)<h(t),
\\
|u(\cdot,t)|_1=1& t>0,
\\
u(\cdot,0)=u_0.
\end{cases}
\end{equation}
The role of $|u|_\iy$ in \eqref{c1} is played here by $h$, and the condition that $\beta_t$ is a probability measure for a.e.\ $t$ is recast by the mass conservation condition.
(For a treatment of an obstacle problem of parabolic type via viscosity solutions, with a given obstacle $h$ and without a global constraint such as $|u(\cdot,t)|_1=1$, see \cite[Theorem 8.6]{el1997reflected}).

Although there are similarities, the tools required are quite different. The proofs in \cite{de-masi-nbbm,de-masi-book, de-masi-excl, de-masi-nbbm2} use the barrier method technique, that is similar to a Trotter scheme which allows one to separate the motion and branching from the removal mechanism. This technique relies crucially on monotonicity properties with respect to mass transport inequalities, which does not seem to have an analogue in the model we study here. We were unable to implement a Trotter scheme for proving uniqueness of solutions to \eqref{c1}. The proofs in \cite{ber19, ber21} are based on the existence of classical solutions to \eqref{n1}.

\subsection{Sketch of the proof}

The proof is based on the compactness--uniqueness approach. 
Showing uniqueness of solutions to \eqref{c1} crucially uses
the continuity of the component $u$ of the solution (specifically, that $u\in\U(\R^d\times\R_+)$). Therefore, relative compactness of the laws of $(\xi^N,\beta^N)$ is insufficient on its own, and must be supplemented with estimates that guarantee that limits of $\xi^N$ have continuous densities w.r.t.\ the Lebesgue measure on $\R^d$.

The first step of the proof, carried out in Section \ref{sec2}, is to take the $N$ limit, with $\eps$ fixed, and show that it can be described by an ODE on $\calS_\eps$, equation \eqref{01}. This equation can be seen as a discrete analogue of equation \eqref{c1} along with the constraint $\beta\in\bB^u$ expressed by \eqref{betau}. In Lemma \ref{lem03}, it is shown that the ODE has at most one solution. Then, in the series of Lemmas \ref{lem05}, \ref{lem06}, \ref{lem04} it is shown that the sequence $(\xi^N,\beta^N)$ is tight. Finally, in Proposition \ref{prop2}, it is shown that all limits are supported on solutions to the ODE \eqref{01}. Hence, denoting in what follows the unique solution to \eqref{01} by $(\bu^{(\eps)},\bLa^{(\eps)})$, the above line achieves the convergence in probability $(\xi^{\eps,N},\beta^{\eps,N})\to(\xi^{(\eps)},\beta^{(\eps)})$ as $N\to\iy$, where
\[
\xi^{(\eps)}_t(dx)=\sum_{y\in\calS_\eps}\bu^{(\eps)}(y,t)\del_y(dx),
\qquad
\beta^{(\eps)}(dx,dt)=\sum_{y\in\calS_\eps}\bLa^{(\eps)}(y,t)\del_y(dx)dt.
\]

The functions $U^{(\eps)}:=\eps^{-d}\bu^{(\eps)}$ defined on $\calS_\eps$ are then considered as prelimit versions for the density of limits of $\xi^{(\eps)}_t$ as $\eps\to0$. The goal of Section \ref{sec3} is to provide estimates on the modulus of continuity of $U^{(\eps)}$ that are uniform in $\eps$ (Proposition \ref{prop1}). The main idea here is to use a coupling of two Markov processes (constructed in equation \eqref{c3}), for which the generator of each is given by $\bar\calL=\Del-b\cdot\grad\cdot$. It is shown in Lemma \ref{lem01} that differences $U^{(\eps)}(x,t)-U^{(\eps)}(y,t)$ can be controlled in terms of this coupling.

The limit as $\eps\to0$ is taken in Section \ref{sec4}, where it is shown, in Proposition \ref{prop4}, that the family $(\xi^{(\eps)},\beta^{(\eps)})$ is relatively compact and limits of its first component have densities that form solutions to \eqref{c1}. Here, the aforementioned estimates on $U^{(\eps)}$ are used to show that these are continuous.

Finally, in Section \ref{sec5}, it is shown that \eqref{c1} has at most one solution. The continuity of the component $u$ of the solution makes it possible to mollify and get approximations in the uniform topology.
We noted above that monotonicity played a role in earlier work in establishing uniqueness for FBP that arise from selection models. The same is true in our approach, although here a different form of monotonicity is essential, namely, the monotonicity of the set-valued operator $u\rightarrow \bB^u$ expressed as follows: $\lan \beta_t-\gamma_t, u(\cdot,t)-v(\cdot,t)\ran \geq 0$ for a.e.\ $t$ provided that $\beta\in\bB^u$ and $\gamma\in\bB^v$.

\subsection{Examples and open problems}
In the discussion that follows, we will refer to the sets $\{x:u(x,t)=|u(\cdot,t)|_\iy\}$ as the {\it argmax sets}. First we show that in some cases explicit solutions can be found. 

\begin{example}[Flat top solution]\label{ex1}
Here we give an explicit stationary solution in dimension 1, in the case where the drift is given by $b(x)=-2\tanh(x)$. Let
\[
w=\log(1+\sqrt{2}),\qquad h=\frac{1}{2(w+\sqrt{2})},
\]
\[
u(x,t)=u(x)=\begin{cases}
h & |x|\le w
\\ \displaystyle
2h\frac{\sinh|x|}{\cosh^2|x|} & |x|>w,
\end{cases}
\]
\[
\beta_t(dx)=h(1+{\rm sech}^2x)1_{[-w,w]}(x) dx, \qquad t\in\R_+.
\]
Then $(u,\beta)$ is a stationary solution of \eqref{c1}, as can be checked by direct calculation. Note that $b$ is bounded Lipschitz, and so our results apply; in particular, the convergence stated in Theorem~\ref{th1} holds provided the initial density is $u_0=u$.
\end{example}

\begin{example}[Sharp and flat top solutions]\label{ex2}
We next consider the drift $b(x)=-a\,{\rm sign}(x)$, $a\ge2$ and show that there are multiple stationary solutions in this case. Some of them have argmax sets in the form of an interval (``flat top''), as in Example \ref{ex1}, while in others the argmax is a singleton (``sharp top'').
Although the discontinuous drift coefficient precludes applying our results, the solutions to \eqref{c1} we construct are well-defined in the sense of Definition~\ref{def1} and provide motivation for some new questions.

First, consider $a>2$. Then for any $w\geq 0$, the following is a stationary solution of \eqref{c1} in the sense of Definition~\ref{def1}:
\begin{equation}\label{10}
	u(x,t) = u(x)= \begin{cases}
	\frac{1}{2(w+a)} & |x|\leq w 
	\\ \displaystyle
	\frac{1}{2(w+a)r} \left[ \lambda_2 e^{-\lambda_1 (|x|-w)} - \lambda_1 e^{-\lambda_2 (|x|-w)} \right] & |x|>w
	\end{cases}
\end{equation}
$$\beta_t(dx)=\frac{a\delta_0(dx)}{w+a} +  \frac{1_{[-w,w]}(x) }{2(w+a)}dx$$
where $r=\sqrt{a^2-4}$ and the decay constants are
$\lambda_1 = \frac{a - r}{2}$, $\lambda_2 = \frac{a + r}{2}$.

%Denote $r=\sqrt{b_0^2-4}$ and
%\[
%\la_i=\frac{-b_0\pm r}{2}
%\qquad
%C_i=\frac{\pm(1-b_0)+r}{2r},
%\qquad i=1,2
%\]
%with $C_1>C_2>0>\la_1>\la_2$. Let
%\[
%u(x,t)=u(x)=C_1e^{\la_1|x|}+C_2e^{\la_2|x|},\qquad x\in\R,
%\]
%\[
%One can check that $(u,\beta)$ is a stationary solution 

%It is remarkable that there are, in fact, infinity many solutions for which  $w=0$ and  $\beta=\delta_{x=0}$.  The solution (\ref{10}) is the only one which has a continuous derivative at $x=0$.

%Another solution is $u_\infty (x):= \frac{e^{-\lambda_2|x|}}{2\lambda_1}$ is a solution with the fastest decay at $\infty$. Indeed, for any $u(0)\in \left[ \frac{1}{2 b_0}, \frac{1}{b_0-\sqrt{b_0^2-4}} \right]$ there exists a unique solution corresponding to $w=0$.

Next, one can show that for any $u(0)\in[ \frac{1}{2 a}, \frac{1}{a-r}]$ there exists a stationary solution $(u,\beta)$ where the argmax set is the singleton $\{0\}$ and $\beta_t(dx)=\del_0(dx)$. One member of this family is \eqref{10} with $w=0$, which is in fact the only one which has a continuous derivative at $x=0$. To mention another member of this family, $u(x)= \frac{e^{-\lambda_2|x|}}{2\lambda_1}$, $\beta_t(dx)=\del_0(dx)$, where $u'$ is discontinuous and $u''$ has an atom.

Finally, for $a=2$, $u(x)=\frac{1}{4}e^{-|x|}(|x|+1)$, $\beta_t(dx)=\del_0(dx)$ is a stationary solution.
%This is a singular case due to the discontinuity of the drift  at $x=0$.
%This example is not covered by our results.  However, it motivates some questions, as follows.
\end{example}

These examples give rise to some questions, as follows.

\begin{problem}
Extend the well-posedness of the PDE and the convergence results beyond (Lipschitz) continuity of $b$.
\end{problem}

\begin{problem}
Is the (Lipschitz) continuity of $b$ sufficient for the argmax sets to have positive Lebesgue measure, and for $\beta$ to be absolutely continuous? More generally, how does the structure of the argmax sets depend on the irregularity of $b$?
\end{problem}

The component $\beta$ of the solution can also be considered as a sub-gradient of the convex operator $F_0:\mathbb{H}^1(\R^d)\to\R$, $F_0(u)= |u|_\infty$.
Focusing for a moment on the case $b=0$, \eqref{c1} can then be considered the gradient flow of the $\la$-convex operator $F(u)=\frac{1}{2}\int(|\grad u|^2-u^2)dx+|u|_\iy$ and one could use the  Br\'ezis-Komura Theorem \cite{brezis2011functional} to prove that \eqref{c1} has at most one solution in the class $\Ll_{2, {\rm loc}}(\R_+,\mathbb{H}^1)$.
Developing a proof of Theorem \ref{th1} using this approach would necessitate estimates to ensure that limit densities lie in $\Ll_{2, {\rm loc}}(\R_+,\mathbb{H}^1)$, that are quite different from those developed here to achieve continuity, but it seems natural to ask about solutions in this space.

\begin{problem} Under what assumptions on $b$ can one obtain existence and uniqueness for \eqref{c1} in $\Ll_{2, {\rm loc}}(\R_+,\mathbb{H}^1)$ and a corresponding convergence result?
\end{problem}

%Another question is about qualitative properties of $\beta$.
%
%\begin{problem}
%Let $(u,\beta)$ solve \eqref{c1}.
%Does $\beta_t$ have density for a.e.\ $t$ under some mild conditions?
%\end{problem}
%
%Next, let $(u,\beta)$ be a solution of \eqref{c1} and denote $h(t)=|u(\cdot,t)|_\iy$. If $O\subset\R^d\times(0,\iy)$ is an open set contained in $\{(x,t):u(x,t)=h(t)\}$ then $\Del u=0$ there, and formally, \eqref{c1} gives $\beta_t(dx)=(-\dot h-h\grad\cdot b(x)+h)1_{\{u(x,t)=h(t)\}}dx$. Moreover, $\beta_t$ is a probability measure for a.e.\ $t$, and consequently $|u(\cdot,t)|_1=1$ for all $t$. This heuristic suggests that \eqref{c1} can be posed as
%\begin{equation}\label{A3}
%\pl_t u=\calL^* u+u+(\dot h+h\grad\cdot b-h)1_{\{u=h\}},
%\qquad h(t)=|u(\cdot,t)|_\iy,
%\qquad |u(\cdot,t)|_1=1
%\end{equation}
%with an initial condition.
%
%\begin{problem} Can one show existence and uniqueness of weak solutions to \eqref{A3} and equivalence to \eqref{c1}?
%\end{problem}

\section{From particle system to ODE}\label{sec2}

In this section it is shown that, for fixed $\eps$
and $N\to\iy$, the empirical measure of the particle system converges to deterministic dynamics given by a system of ODE indexed by $\calS_\eps$.
Throughout the section, $\eps$ is fixed.
Moreover, Assumptions \ref{assn1} and \ref{assn2} apply throughout the remainder of this article. Recall that by Assumption \ref{assn2},
we have $\xi^{\eps,N}_0\to\xi^{(\eps)}_0$ in probability in $\calM(\R^d)$ as $N\to\iy$.

First we define the formal adjoint of $\calL_\eps$ and introduce some related notation.
For $1\le i\le 2d$, let $i^*=i+d$ if $i\le d$
and $i^*=i-d$ if $i\ge d+1$. With this, $k_{i^*}=-k_i$.
Let
\[
\rho_{\eps,i}(x)=r_{\eps,i^*}(x+\eps k_i).
\]
Denote $\bar r_{\eps}(x)=\sum_{i\le 2d}r_{\eps,i}(x)$
and $\bar\rho_\eps(x)=\sum_{i\le2d}\rho_{\eps,i}(x)$.
Let
\begin{align*}
\calL^*_\eps f(x) &= \sum_{i\le2d}
[r_{\eps,i}(x-\eps k_i)f(x-\eps k_i)-r_{\eps,i}(x)f(x)]
\\
&=\sum_{i\le 2d}\rho_{\eps,i}(x)f(x+\eps k_i)-\bar r_\eps(x)f(x),
\qquad\quad x\in\calS_\eps.
\end{align*}
The duality relation holds,
\begin{equation}\label{q4}
\lan \calL_\eps f,g\ran=\lan f,\calL_\eps^*g\ran,
\end{equation}
for $f$ and $g$ defined on $\calS_\eps$, and, say, one of them is compactly supported, and $\lan\cdot,\cdot\ran$ denoting the usual inner product on $\calS_\eps$.

The ODE on $\calS_\eps$ is an equation for $(\bu,\bLa)$, which, given an initial condition $\bu_0\in\Ll_1(\calS_\eps,\R_+)$, reads
\begin{equation}\label{01}
\begin{split}
& \bu(x,t)=\bu_0(x)+\int_0^t(\calL^*_\eps \bu+\bu-\bLa)(x,s)ds,
\quad(x,t)\in\calS_\eps\times\R_+\\
&\bLa(x,t)\ge0, \quad
\sum_x\bLa(x,t)=1,
 \hspace{7.2em} x\in\calS,\ \text{a.e. } t\in\R_+,
\\
&\sum_x\int_0^\iy(|\bu(\cdot,t)|_\iy-\bu(x,t))\bLa(x,t)dt=0.
\end{split}
\end{equation}
\begin{definition}[Solution to \eqref{01}]\label{def2}
A solution to \eqref{01} is a pair $(\bu,\bLa)$ where $\bu$ is a mapping $\calS_\eps\times\R_+\to\R_+$, $\bu(x,\cdot)\in C(\R_+,\R_+)$ for all $x$, $|\bu(\cdot,t)|_1=1$ for all $t$, $\bLa$ is a measurable mapping $\calS_\eps\times\R_+\to\R_+$, and \eqref{01} holds.
\end{definition}

\begin{lemma}\label{lem03}
i.
If $(\bu,\bLa)$ and $(\bar\bu,\bar\bLa)$ are two solutions to \eqref{01} for an initial condition $\bu_0\in\Ll_1(\calS_\eps,\R_+)$ then $\bu=\bar\bu$ and, for every $x$ and a.e.\ $t$, $\bLa(x,t)=\bar\bLa(x,t)$.
\\
ii. Suppose there exists a solution $(\bu,\bLa)$ to \eqref{01}. Then for $0\le s<t$,
\begin{equation}\label{t0}
\bu(y,t)=\sum_{x\in\calS_\eps}\bu(x,s)\frs_{t-s}(x,y)-\int_s^t\sum_{x\in\calS_\eps}\bLa(x,\tau)\frs_{t-\tau}(x,y)d\tau,
\end{equation}
where $\frs_t(x,y)=\frs^{(\eps)}_t(x,y)=e^t\frp_t(x,t)$, and $\frp=\frp^{(\eps)}$ is the fundamental solution of $\pl_t\bu=\calL^*_\eps\bu$ on $\calS_\eps$.
\end{lemma}

\proof
i.
Denote $\bw=\bu-\bar\bu$. Clearly, $|\bu(\cdot,t)|_2\vee|\bar\bu(\cdot,t)|_2\le1$. Moreover,
\begin{align*}
|\bw(\cdot,t)|^2_2 &= |\bw(\cdot,0)|^2_2+2\sum_x\int_0^t\bw(\calL^*_\eps\bw+\bw-\bLa+\bar\bLa)(x,s)ds.
\end{align*}
By Assumption \ref{assn1}, $r_\eps$ is bounded (recall $\eps$ is fixed). Hence
$|\calL^*_\eps f(x)|\le c\sum_{i\le2d}|f(x+\eps k_i)|+c|f(x)|$, and
\[
\sum_x|f(x)\calL^*_\eps f(x)|\le c|f|^2_2.
\]
This shows that
\[
|\bw(\cdot,t)|^2_2\le c\int_0^t|\bw(\cdot,s)|^2_2ds+\int_0^t\sum_x(\bu-\bar\bu)(-\bLa+\bar\bLa)(x,s)ds.
\]
The last integral is $\le0$ because for all $s$,
\[
\sum_x\bu(x,s)\bar\bLa(x,s)
\le \sum_x|\bu(\cdot,s)|_\iy\bar\bLa(x,s) = |\bu(\cdot,s)|_\iy
\]
whereas for a.e.\ $s$,
\[
\sum_x\bu(x,s)\bLa(x,s)
= \sum_x|\bu(\cdot,s)|_\iy\bLa(x,s) = |\bu(\cdot,s)|_\iy,
\]
and a similar statement holds for $\bar\bu$. Hence by Gronwall's lemma, $\bw=0$.

The a.e.\ equality $\bLa=\bar\bLa$ follows because, by \eqref{01}, for a.e.\ $t$, $\bLa(x,t)=(\dot\bu-\calL_\eps^*\bu-\bu)(x,t)$.

ii. It is easy to see that $\frs^{(\eps)}$ is the fundamental solution of $\pl_t\bu=\calL^*_\eps\bu+\bu$ on $\calS_\eps$. Thus \eqref{t0} follows from the first line in \eqref{01} by Duhamel's principle.
\qed

In the remainder of this section we show that $(\xi^N,\beta^N)$ is a tight sequence and its limit can be characterized in terms of the unique solution to \eqref{01}.
A useful tool here will be a coupling of $\xi^N$ of \eqref{20a} with a process in which removals do not take place, namely
\begin{align}\notag
\bar\xi^N_t=\xi^N_0
&+\frac{1}{N}\int_{\calS_\eps\times[2d]\times\R_+\times[0,t]}1_{[0,Nr_i(x)\bar\xi^N_{s-}(x)]}(\theta)
(\del_{x+\eps k_i}-\del_x)
\Pi(dx,di,d\theta,ds)
\\ \label{20b}
&+\frac{1}{N}\int_{\calS_\eps\times\R_+\times[0,t]} 1_{[0,N\bar\xi^N_{s-}(x)]}(\theta)\del_x\,\Gamma(dx,d\theta,ds).
\end{align}
This is the normalized configuration measure for mutually independent branching random walks. It is a standard fact that $\bar\xi^N$ dominates $\xi^N$, but for completeness, we provide a proof.

\begin{lemma}\label{lem05} One has
$\xi^N_t(x)\le\bar\xi^N_t(x)$ for all $t$ and $x$ a.s.
\end{lemma}
\proof
Let
\[
\tau=\inf\{t:\bar\xi^N_t(x)<\xi^N_t(x) \text{ for some } x\}.
\]
On an event of full measure, $\Pi$ and $\Gamma$ are locally finite, and charge at most one point mass to each set of the form $\calS_\eps\times[2d]\times\R_+\times\{t\}$ or $\calS_\eps\times\R_+\times\{t\}$, respectively.
Consider the intersection of this event with $\{\tau<\iy\}$. On this event, $\xi^N_t(x)\le\bar\xi^N_t(x)$ for all $t<\tau$ and all $x$, and moreover, there exists $y$ such that $\bar\xi^N_\tau(y)<\xi^N_\tau(y)$. Note that $\bar\xi^N_{\tau-}(y)=\xi^N_{\tau-}(y)$ must hold.
This can only occur if (the cases below correspond to the first three terms in \eqref{20a}; the fourth term does not have a counterpart of in \eqref{20b}):

- $\Pi$ has, for some $i$ and some $x$, an atom in $\{x\}\times\{i\}\times[0,Nr_i(x)\xi^N_{\tau-}(x)]\times\{\tau\}$, but it does not have one in $\{x\}\times\{i\}\times[0,Nr_i(x)\bar\xi^N_{\tau-}(x)]\times\{\tau\}$. This is impossible because $\xi^N_{\tau-}(x)\le\bar\xi^N_{\tau-}(x)$.

- $\Pi$ has, for some $i$, an atom in $\{y\}\times\{i\}\times[0,Nr_i(y)\bar\xi^N_{\tau-}(y)]\times\{\tau\}$, but it does not have one in $\{y\}\times\{i\}\times[0,Nr_i(y)\xi^N_{\tau-}(y)]\times\{\tau\}$. This is impossible because $\bar\xi^N_{\tau-}(y)=\bar\xi^N_{\tau-}(y)$.

- $\Gamma$ has an atom in $\{y\}\times[0,Nr_i(y)\xi^N_{\tau-}(y)]\times\{\tau\}$, but it does not have one in $\{y\}\times[0,Nr_i(y)\bar\xi^N_{\tau-}(y)]\times\{\tau\}$. This is impossible because $\xi^N_{\tau-}(y)\le\bar\xi^N_{\tau-}(y)$.

Thus $\PP(\tau<\iy)=0$ and the result follows.
\qed

Using this we can provide an upper bound on $\xi^N$ as follows.
Let $\calB_r=\{x\in\R^d:|x|<r\}$.

\begin{lemma}\label{lem06}
Given $T>0$, $\del>0$ there exist $r>0$ such that
\[
\sup_N\sup_{t\in[0,T]}\E[\xi^N_t(\calB_r^c)] \le \del.
\]
\end{lemma}

\proof
In view of Lemma \ref{lem05}, it suffices to prove the estimate for $\bar\xi^N$. Now, $\bar\xi^N$ is nothing but the normalized configuration measure of independent branching RWs starting $\xi^N_0$.
To be more precise, let $\bar X^{N,i}_t$, $i\in[N]$ be a collection of Branching random walks on $\calS_\eps$ that are mutually independent conditioned on the initial configuration $\sig\{\bar X^{N,i}_0\}$, with branching at rate $1$ and generator $\calL_\eps$. Let the root particle of $\bar X^{N,i}$ be labeled as $(i,0)$ and the descendants as $(i,j)$, $j\ge1$. Let $\tau^{N,i,j}$ be the birth time of particle $\bar X^{N,i,j}$. Let
\[
\tilde\xi^N_t=\frac{1}{N}\sum_{i,j}1_{\{t\ge\tau^{N,i,j}\}}\del_{\bar X^{N,i,j}_t}
\]
with the initial condition $\tilde\xi^N_0=\xi^N_0$. Then the process $\tilde\xi^N$ is equal in law to $\bar\xi^N$. Hence by the many-to-one lemma \cite{harris2017},
\begin{align*}
\E[\bar\xi^N_t(\calB_r^c)]&=\frac{1}{N}\sum_{i=1}^N\E\Big[\sum_j1_{\{t\ge\tau^{N,i,j}\}}1_{\calB_r^c}(\bar X^{N,i,j}_t) \Big]
= \frac{e^t}{N}\sum_{i=1}^N \E[1_{\calB_r^c}(\bar X^{N,i,0}_t)]
\\
&\le
e^t\E[\xi^N_0(\calB_{r/2}^c)]+\frac{e^t}{N}\sum_{i=1}^N\PP\Big(|\bar X^{N,i,0}_t-\bar X^{N,i,0}_0|\ge \frac{r}{2}\Big).
\end{align*}
The above expression is bounded, uniformly in $(N,t)\in\N\times[0,T]$, by $\del$, provided that $r$ is sufficiently large. Indeed, the first term is bounded by $\del/2$ provided $r$ is large, owing to the tightness of $\xi^N_0$, $N\in\N$ (for each $\eps$) in $\calM(\R^d)$, implied by Assumption \ref{assn2}. For the second term, the same is true because the number of jumps $\bar X^{N,i,0}$ makes during $[0,t]$ is dominated by a Poisson$(cT)$ where $c$ does not depend on $i$ (but may depend on $\eps$). The result follows.
\qed

Next, a calculation.
Consider the marginals $\Pi_{x,i}$, $\Gamma_x$ defined as
\begin{align*}
\int_{\R_+^2}g(\theta,t)\Pi_{x,i}(d\theta,dt)&=\int_{\calS_\eps\times[2d]\times\R_+^2} g(\theta,t)1_{\{y=x,j=i\}}\Pi(dy,dj,d\theta,dt),
\qquad x\in\calS_\eps, \ i\in[2k],
\\
\int_{\R_+^2}g(\theta,t)\Gamma_x(d\theta,dt)&=\int_{\calS_\eps\times\R_+^2} g(\theta,t)1_{\{y=x\}}\Pi(dy,d\theta,dt),
\hspace{6.4em} x\in\calS_\eps.
\end{align*}
Then by \eqref{20a},
\begin{align}\label{21a}
\notag
u^N(x,t) &= u^N_0(x)-\frac{1}{N}\sum_i\int_{\R_+\times[0,t]}1_{[0,Nr_i(x)u^N(x,s-)]}(\theta)\Pi_{x,i}(d\theta,ds)
\\ \notag
&\qquad +\frac{1}{N}\sum_i\int_{\R_+\times[0,t]}1_{[0,Nr_i(x-\eps k_i)u^N(x-\eps k_i,s-)]}(\theta)\Pi_{x-\eps k_i,i}(d\theta,ds)
\\ \notag
&\qquad + \frac{1}{N}\int_{\R_+\times[0,t]}1_{[0,Nu^N(x,s-)]}(\theta) \Gamma_x(d\theta,ds)
\\
&\qquad -\frac{1}{N}\sum_y \int_{\R_+\times[0,t]}1_{[0,Nu^N(y,s-)]}(\theta)1_{\{x=\argmax^*u^N(\cdot,s-)\}} \Gamma_y(d\theta,ds).
\end{align}
Construct the signed measures
\[
\hat\Pi_{x,i}(d\theta,dt)=\Pi_{x,i}(d\theta,dt)-d\theta\,dt,
\qquad
\hat\Gamma_x(d\theta,dt)=\Gamma_x(d\theta,dt)-d\theta\,dt.
\]
Then
\begin{align}\label{25}
\notag
u^N(x,t)&=u^N_0(x)+\int_0^t\sum_i [r_i(x-\eps k_i)u^N(x-\eps k_i,s)-r_i(x)u^N(x,s)]ds
\\ \notag
&\hspace{4em}
+\int_0^tu^N(x,s)ds
-\sum_y\int_0^t u^N(y,s)
1_{\{x=\argmax^*u^N(\cdot,s)\}}ds+M^N_x(t)
\\ \notag
&=u^N_0(x)+\int_0^t [(\calL^*_\eps u^N+u^N)(x,s)
-1_{\{x=\argmax^*u^N(\cdot,s)\}}]ds+M^N_x(t)
\\
&=u^N_0(x)+\int_0^t [(\calL^*_\eps u^N+u^N)(x,s)]ds
-\beta^N(\{x\}\times[0,t])+M^N_x(t),
\end{align}
where $M^N_x$ is obtained by replacing $\Pi$ and $\Gamma$ by $\hat\Pi$ and $\hat\Gamma$ on the right-hand side of \eqref{21a}. This is a martingale for each $x$.
To provide a uniform (in $x$) estimate on its quadratic variation, fix $R$ and $T$. Let $c_1=c_1(\eps)$ be an upper bound on $r_i(x)$ for $x\in \calB_{R+1}$. For $(x,t)\in \calB_R\times[0,T]$, using $\sum_yu^N(y,t)=1$,
\begin{equation}\label{26}
\sum_{x\in\calS_\eps}\E\{[M^N_x]_T\} \le N^{-2}(2dNc_1+2dNc_1+N+N)T
=c_2TN^{-1}.
\end{equation}

\begin{lemma}\label{lem04}
$\xi^N$ is $C$-tight in $D(\R_+,\calM(\calS_\eps))$
and $\beta^N$ is tight in $\calM^{(1)}(\calS_\eps\times\R_+)$.
\end{lemma}

\proof
A sufficient condition for tightness of $\beta^N$ as a sequence of $\calM^{(1)}(\calS_\eps\times\R_+)$-valued random variables, is that for every $\del>0$ and $T$ there is $r$ such that
\begin{equation}\label{23a}
\limsup_N\PP(\beta^N(\calB_r^c\times[0,T])>\del)<\del.
\end{equation}
By~\eqref{25},  Fubini, and the fact that $\calL^*_\eps$ is nearest neighbor bounded operator, we have
\begin{align}
\notag
\E\left[ \beta^N(\{\calB_r^c\times[0,T]\})\right] &\leq \E [u^N_0(\calB_r^c)]+\E\left[ \int_0^T \sum_{x\in \calB_r^c}   [(\calL^*_\eps u^N+u^N)(x,s)]ds\right].
\\
\notag
& 
\leq \E [u^N_0(\calB_r^c)]  + C_{\eps,T }\int_0^T \sum_{x\in \calB_{r/2}^c} \E\left[ u^N(x,s)\right]ds.
\end{align}
Now, by Lemma~\ref{lem06} we easily obtain that for every $\del'>0$ there exists $r>0$ such that for all $N\geq 1$,
\begin{align}
\notag
\E\left[ \beta^N(\{\calB_r^c\times[0,t]\})\right] \leq \delta'/2+\delta'/2=\delta',
\end{align}
and \eqref{23a} follows.

As for $C$-tightness of the $\calM_1(\calS_\eps)$-valued processes $\xi^N$, it suffices to prove
that for every $T$ and $\del>0$ there is $r$ such that
\begin{equation}\label{26a}
\limsup_N\sup_{t\in[0,T]}\PP(\xi^N_t(\calB_r^c)>\del)<\del,
\end{equation}
and for every $T$, $x\in\calS_\eps$ and $\del>0$ there is $\eta>0$ such that
\begin{equation}\label{27a}
\limsup_N\PP(\sup_{x\in\calS_\eps}\om_T(\xi^N_\cdot(x),\eta)>\del)<\del.
\end{equation}
Here, $\om_T(f,\eta)=\sup\{|f(s)-f(t)|:s,t\in[0,T], |s-t|\le\eta\}$.
Note that \eqref{26a} follows from Lemma \ref{lem06}.
To prove \eqref{27a}, let $0\le s\le t\le T$. Because the rates $r_\eps$ are bounded, we have by \eqref{25} that
\[
|u^N(x,t)-u^N(x,s)|\le c(t-s)+|M^N_x(t)-M^N_x(s)|.
\]
Hence the estimate \eqref{26} on the quadratic variation of $M^N_x$ implies \eqref{27a}.
\qed

We can now state the main result of this section, regarding the convergence as $N\to\iy$ for fixed $\eps$.
\begin{proposition}\label{prop2}
There exists $(\xi,\beta)\in C(\R_+,\calM(\calS_\eps))\times\calM^{(1)}(\calS_\eps\times\R_+)$ such that $(\xi^N,\beta^N)\to(\xi,\beta)$ in $D(\R_+,\calM(\calS_\eps))\times\calM^{(1)}(\calS_\eps\times\R_+)$ in probability.
Moreover, define $(\bu,\bLa)$ via
\begin{equation}\label{32a}
\xi_t(\{x\})=\bu(x,t),
\qquad
\beta(dx,dt)=\sum_{y\in\calS_\eps}\bLa(y,t)\del_y(dx)dt.
\end{equation}
Then $(\bu,\bLa)$ is a solution to \eqref{01} (which is unique according to Lemma \ref{lem03}).
\end{proposition}

\proof
In view of the tightness stated in Lemma \ref{lem04} and the uniqueness stated in Lemma \ref{lem03}, it suffices to show that whenever $(\xi^N,\beta^N)\to(\xi,\beta)$ along a subsequence, $(\bu,\bLa)$ defined via \eqref{32a} solves \eqref{01} in the sense of Definition \ref{def2}.

To this end, note that the continuity of $\bu(x,\cdot)$ for every $x$ follows from $C$-tightness of $\xi^N$.
Next, recalling the estimate \eqref{26} and invoking Skorohod's representation, we assume without loss of generality that, along the subsequence, $(\xi^N,\beta^N)\to(\xi,\beta)$ a.s.\ and $M^N_x\to0$ for all $x$, a.s. Now take the limit in \eqref{25}. Then, for every $x$ and $t$, $u^N(x,t)\to\bu(x,t)$, and by bounded convergence, the integral in \eqref{25} converges to $\int_0^t(\calL^*_\eps \bu+\bu)(x,s)ds$. We obtain that the first line in \eqref{01} is satisfied.

To show that the last line in \eqref{01} holds, note that, by construction, for every $x$,
\begin{equation}\label{w30}
\int_0^\iy(|u^N(\cdot,s)|_\iy-u^N(x,s))\beta^N_x(ds)=0,
\end{equation}
where $\beta^N_x$ is the marginal at $x$.
Arguing by contradiction, suppose there exists an $x$ such that
$\int_0^\iy(|\bu(\cdot,s)|_\iy-\bu(x,s))\bLa(x,s)ds>0$. By continuity of $\bu$ in the second variable, there exists an interval $(t_1,t_2)$ on which $|\bu(\cdot,s)|_\iy-\bu(x,s)>\del>0$ and $\int_{t_1}^{t_2}\bLa(x,s)ds>\del$. Because $\xi^N\to\xi$, with $\xi$ a continuous path, we have u.o.c.\ convergence of $u^N$, hence for all large $N$,
\[
|u^N(\cdot,s)|_\iy-u^N(x,s) > \del/2, \qquad s\in(t_1,t_2).
\]
Moreover, $\beta^N_x(dt)\to\bLa(x,t)dt$ as measures, hence $\beta^N_x(t_1,t_2)>\del/2$ for all large $N$.
Hence
\[
\int_{t_1}^{t_2}(|u^N(\cdot,s)|_\iy-u^N(x,s))\beta^N_x(ds)>\del^2/4,
\]
contradicting \eqref{w30}. The result follows.
\qed

\section{Estimates on modulus of continuity}\label{sec3}

The goal of this section is to provide an estimate, that is uniform in $\eps$, on the modulus of continuity of $U^{(\eps)}:=\eps^{-d}\bu^{(\eps)}$, where $\bu^{(\eps)}$ is the first component of the solution of \eqref{01}. This will enable us to show that limits of $\xi^{N,\eps}$ (as $N\to\iy$ followed by $\eps\to0$) have continuous densities.

For $f:A\to\R$, $A\subset\R^d$, let
\[
\om(f,\eta)=\sup\{|f(x)-f(y)|: x,y\in A,\,|x-y|\le\eta\}.
\]
For $f:\calS_\eps\times[0,\iy)\to\R$,
\[
\pl^{(\eps)}_if(x,t) = f(x+\eps k_i,t)-f(x,t).
\]
For the most part, we will suppress the dependence on
$\eps$ and write $\pl_i$ for $\pl_i^{(\eps)}$. Note that
$\calL^*_\eps$ is not the generator of a Markov process. However, if we let
\begin{equation}\label{03}
\bar\calL_\eps f(x)=\sum_{i\le2d}\rho_{\eps,i}(x)\pl_if(x), \qquad x\in\calS_\eps,
\end{equation}
then it is a generator, and it is related to $\calL^*_\eps$
via
\[
\calL^*_\eps f(x)=\bar\calL_\eps f(x)+(\bar\rho_\eps(x)-\bar r_\eps(x))f(x).
\]
Denote $h_\eps(x)=\bar r_\eps(x)-\bar\rho_\eps(x)$.
\begin{lemma}\label{lem0}
There exists a constant $C_2$ such that
\[
|h_\eps|_\iy\le C_2,
\qquad
\sup_{x,y\in\calS_\eps,x\ne y}\frac{|h_\eps(x)-h_\eps(y)|}{|x-y|}\le C_2,\qquad \eps\in(0,1).
\]
\end{lemma}
\proof Suppressing $\eps$,
\begin{align*}
h(x)&=\bar r(x)-\bar\rho(x)=\sum_{i\le2d}(r_i(x)-r_{i^*}(x+\eps k_i))=\sum_{i\le2d}(r_i(x)-r_i(x-\eps k_i))
\\
&=\eps^{-1}\sum_{i\le2d}(q_i(x)-q_i(x-\eps k_i)),
\end{align*}
and the bound on $|h|$ follows by the first part of \eqref{c01}. Similarly,
\[
h(x)-h(y)=\eps^{-1}\sum_{i\le2d}[(q_i(x)-q_i(x-\eps k_i))
-(q_i(y)-q_i(y-\eps k_i))],
\]
and the second assertion
of the lemma follows from the second part of \eqref{c01}.
\qed

Let $(\bar X_t,\bar Y_t)$ be a coupling of Markov processes
on $\calS_\eps$ each governed by the generator $\bar\calL_\eps$, constructed as follows.
Let $\bar\Pi_i$, $1\le i\le 2d$ be
mutually independent PRMs on
$\R_+^2$, with the Lebesgue
measure on $\R_+^2$ as the intensity measure. We use the same collection $\bar\Pi_i$
to drive both processes $X$ and $Y$, namely
\begin{equation}\label{c3}
\begin{split}
\bar X_t&=\bar x+\sum_{i\le 2d}\eps k_i\int_{[0,t]\times\R_+}
1_{[0,\rho_i(\bar X_{s-})]}(\theta)\bar\Pi_i(ds,d\theta),
\\
\bar Y_t&=\bar y+\sum_{i\le 2d}\eps k_i\int_{[0,t]\times\R_+}
1_{[0,\rho_i(\bar Y_{s-})]}(\theta)\bar\Pi_i(ds,d\theta).
\end{split}
\end{equation}

For $\eps>0$ and with $\bu^{(\eps)}_0$ as in Assumption \ref{assn2}, let $(\bu^{(\eps)},\bLa^{(\eps)})$ be the solution of \eqref{01} corresponding to $\bu^{(\eps)}_0$, and let $U^{(\eps)}:=\eps^{-d}\bu^{(\eps)}$. By Assumption \ref{assn2}, $U^{(\eps)}$ are uniformly bounded, and
\begin{equation}\label{12w}
\lim_{\eta\to0}\limsup_{\eps\to0}\om(U^{(\eps)}_0,\eta)=0.
\end{equation}
The goal of this section is to show that the modulus of continuity of $U^{(\eps)}(\cdot,t)$ can be controlled uniformly in $\eps$, as follows.

\begin{proposition}\label{prop1}
For every $T$,
\begin{equation}\label{05}
\lim_{\eta\to0}\limsup_{\eps\to0} \sup_{t\in(0,T]} \om(U^{(\eps)}(\cdot,t),\eta)=0.
\end{equation}
\end{proposition}

The main step toward proving Proposition \ref{prop1} is the following estimate.
Recall the constants $C_1$, $C_2$ from Assumption \ref{assn1} and Lemma \ref{lem0}.
\begin{lemma}\label{lem01}
(a) One has $|U^{(\eps)}(\cdot,t)|_\iy\le |U^{(\eps)}_0|_\iy \, e^{(1+C_2)t}$.
\\
(b)
Let $\bar x,\bar y\in\calS^{(\eps)}$ and let
$(\bar X_t,\bar Y_t)$ be the coupling \eqref{c3}. Let $\tau$ be a stopping time w.r.t.\ the filtration generated by $(\bar X_t,\bar Y_t)$. Then for every $T>0$ and $0\le t\le T$,
\[
(U^{(\eps)}(\bar x,t)-U^{(\eps)}(\bar y,t))^2
\le c_T\Big[\E(U^{(\eps)}_0(\bar X_{t\w\tau})-U^{(\eps)}_0(\bar Y_{t\w\tau}))^2+\E\int_0^t|\bar X_{s\w\tau}-\bar Y_{s\w\tau}|^2ds\Big],
\]
where $c_T$ depends only on $C_1$, $C_2$ and $T$.
\end{lemma}

\proof
(a)
Suppress $\eps$ from the notation of $\bu,\bLa,\rho$, etc.
Fix $T>0$. For $t\in[0,T]$, let $\bv(x,t)=\bu(x,T-t)$.
Then
\begin{equation}\label{02}
\dot \bv(x,t)=-\calL^*_\eps \bv(x,t)-\bv(x,t)+\bLa(x,T-t).
\end{equation}
We have
\begin{align*}
\bv(\bar X_t,t)&=\bv(\bar x,0)+\int_0^t\dot \bv(\bar X_s,s)ds+
\sum_{i\le2d}\int_{[0,t]\times\R_+}\bar\Pi_i(ds,d\theta)
\pl_i\bv(\bar X_{s-},s)1_{[0,\rho_i(\bar X_{s-})]}(\theta)
\\
&=\bv(\bar x,0)+\int_0^t\dot \bv(\bar X_s,s)ds+\sum_{i\le2d}\int_0^t\rho_i(\bar X_s)\pl_i\bv(\bar X_s,s)ds+m_t,
\end{align*}
where $m$ is a martingale.
Now,
\[
\dot \bv+\sum_i\rho_i\pl_i\bv=\dot \bv+\bar\calL \bv=-\bv+\bLa-h\bv
\ge -(1+C_2)\bv.
\]
This gives $\E \bv(\bar X_T,T)\ge \bv(\bar x,0)-(1+C_2)\int_0^T
\E \bv(\bar X_s,s)ds$. Thus, with $a_t=|\bu(\cdot,t)|_\iy$, we have
$a_T\le a_0+(1+C_2)\int_0^Ta_s ds$. This holds for every $T$, and therefore $a_T\le a_0 e^{(1+C_2)T}$.
The result follows.

(b)
For $f:\calS_\eps\times\calS_\eps\times[0,\iy)\to\R$, denote
\begin{align*}
\pl_{1,i}f(x,y,t) &= f(x+\eps k_i,y,t)-f(x,y,t),
\\
\pl_{2,i}f(x,y,t) &= f(x,y+\eps k_i,t)-f(x,y,t),
\\
\pl_{12,i}f(x,y,t) &= f(x+\eps k_i,y+\eps k_i,t)-f(x,y,t).
\end{align*}
By convention, the interval
$(a,b]$ is the empty set if $a\ge b$.
Then for $a,b\ge0$, the interval $[0,a\vee b]$ can be partitioned into $[0,a\w b]$, $(a,b]$ and $(b,a]$, where
at least one of these intervals is empty.
Noting that points
$(s,\theta)$ of $\bar\Pi_i$ with
$\theta\le \rho_i(\bar X_{s-})\w \rho_i(\bar Y_{s-})$ give rise to a jump
of size $\eps k_i$ in {\it both} $\bar X$ and $\bar Y$, we can write, for $f(\cdot,\cdot,\cdot):\calS_\eps\times\calS_\eps\times[0,\iy)$
that is $C^1$ in the last variable,
\begin{align*}
f(\bar X_t,\bar Y_t,t)&=f(\bar x,\bar y,0) + \int_0^t\dot f(\bar X_s,\bar Y_s,s)ds
+\sum_{s\in(0,t]}(f(\bar X_s,\bar Y_s,s)-f(\bar X_{s-},\bar Y_{s-},s))
\\
&=f(\bar x,\bar y,0)+ \int_0^t\dot f(\bar X_s,\bar Y_s,s)ds
\\
&\quad+\sum_{i\le2d}\int_{[0,t]\times\R_+}
\bar\Pi_i(ds,d\theta)
\Big\{\pl_{12,i}f(\bar X_{s-},\bar Y_{s-},s)1_{[0,\rho_i(\bar X_{s-})\w \rho_i(\bar Y_{s-})]}(\theta)
\\
&\hspace{2em}
+\pl_{1,i}f(\bar X_{s-},\bar Y_{s-},s)1_{(\rho_i(\bar Y_{s-}),\rho_i(\bar X_{s-})]}(\theta)
+\pl_{2,i}f(\bar X_{s-},\bar Y_{s-},s)1_{(\rho_i(\bar X_{s-}),\rho_i(\bar Y_{s-})]}(\theta) \Big\}.
\end{align*}
Thus
\begin{equation}\label{04}
f(\bar X_t,\bar Y_t,t)=f(\bar x,\bar y,0)
+\int_0^t F(\bar X_s,\bar Y_s,s)ds +m_t,
\end{equation}
where $m$ is a martingale and
\begin{align*}
F(x,y,t)&=\dot f(x,y,t)
+\sum_{i\le2d}\Big\{\pl_{12,i}f(x,y,t)(\rho_i(x)\w \rho_i(y))
\\
&\hspace{2em}
+\pl_{1,i}f(x,y,t)(\rho_i(x)-\rho_i(y))^+
+\pl_{2,i}f(x,y,t)(\rho_i(y)-\rho_i(x))^+ \Big\}.
\end{align*}

We now fix $T_0>0$ and consider $0<T\le T_0$.
Define $\bv$ in terms of $T$ as in part (a) of the proof.
We will use equation \eqref{04} with
\[
f(x,y,t)=(\bv(x,t)-\bv(y,t))^2, \qquad 0\le t\le T.
\]
Using $a^2-b^2\ge 2b(a-b)$, and suppressing
$t$ from the notation,
\begin{align*}
\pl_{12,i}f(x,y) &\ge2(\bv(x)-\bv(y))(\bv(x+\eps k_i)-\bv(y+\eps k_i)-\bv(x)+\bv(y))
\\
&=2(\bv(x)-\bv(y))(\pl_i\bv(x)-\pl_i\bv(y))
\\
\pl_{1,i}f(x,y) &\ge2(\bv(x)-\bv(y))(\bv(x+\eps k_i)-\bv(y)-\bv(x)+\bv(y))
\\
&=2(\bv(x)-\bv(y))\pl_i\bv(x)
\\
\pl_{2,i}f(x,y) &\ge2(\bv(x)-\bv(y))(\bv(x)-\bv(y+\eps k_i)-\bv(x)+\bv(y))
\\
&=-2(\bv(x)-\bv(y))\pl_i\bv(y).
\end{align*}
Thus (suppressing $t$),
\begin{align*}
& F(x,y)
\ge 2(\bv(x)-\bv(y))\Big\{\dot \bv(x)-\dot \bv(y)
+\sum_{i\le 2d}
\Big[
(\pl_i \bv(x)-\pl_i\bv(y))(\rho_i(x)\w \rho_i(y))
\\
&\hspace{18em}
+\pl_i\bv(x)(\rho_i(x)-\rho_i(y))^+
-\pl_i\bv(y)(\rho_i(y)-\rho_i(x))^+
\Big] \Big\}
\\
&=2(\bv(x)-\bv(y))\Big\{\dot \bv(x)-\dot \bv(y)
+\sum_{i\le 2d} (\pl_i \bv(x)\rho_i(x)-\pl_i\bv(y)\rho_i(y))
\Big\}
\\
&=2(\bv(x)-\bv(y))\Big\{\dot \bv(x)-\dot \bv(y)
+\bar\calL_\eps \bv(x)-\bar\calL_\eps \bv(y)
\Big\},
\end{align*}
where \eqref{03} is used on the last line. Now, by \eqref{02},
\[
(\dot \bv+\bar\calL_\eps \bv)(x,t)=\bLa(x,T-t)
+(h \bv-\bv)(x,t).
\]
Thus
\begin{align*}
F(x,y,t)
&\ge2(\bv(x,t)-\bv(y,t))(\bLa(x,T-t)-\bLa(y,T-t)+(h \bv-\bv)(x,t)
-(h \bv-\bv)(y,t)).
\end{align*}
If $\bLa(x,T-t)-\bLa(y,T-t)>0$ then $\bLa(x,T-t)>0$,
implying that $\bu(x,T-t)-\bu(y,T-t)\ge0$, by \eqref{01},
that is, $\bv(x,t)-\bv(y,t)\ge0$.
Similarity, negativity of the former
implies nonpositivity of the latter. Hence
\begin{align*}
F(x,y,t) &\ge 2(\bv(x)-\bv(y))((h \bv-\bv)(x)-(h \bv-\bv)(y))
\\
&= -2(\bv(x)-\bv(y))^2+2(\bv(x)-\bv(y))[h (x)(\bv(x)-\bv(y))+(h(x)-h(y))\bv(y)]
\\
&\ge -2(1+|h(x)|)(\bv(x)-\bv(y))^2-2|\bv(x)-\bv(y)||h(x)-h(y)| \bv(y)
\\
&\ge -2(1+C_2)(\bv(x)-\bv(y))^2-2C_2|\bv(x)-\bv(y)||x-y|\bv(y).
\end{align*}
With $V=V^{(\eps)}:=\eps^{-d}\bv^{(\eps)}$,
\begin{align*}
\eps^{-2d}F(x,y,t)&\ge-2(1+C_2)(V(x)-V(y))^2-2C_2|V(x)-V(y)||x-y|V(y)
\\
&\ge-c_1(T)(V(x)-V(y))^2-c_2(T)(x-y)^2,
\end{align*}
where $c_1(T)=2(1+C_2)+C_2C_1e^{(1+C_2)T}$, $c_2(T)=C_2C_1e^{(1+C_2)T}$.
Let $c=c(T_0)=c_1(T_0)\vee c_2(T_0)$ and recall that $0\le t\le T\le T_0$.
Then using the above inequality in \eqref{04}, substituting $T\w\tau$ for $t$, we obtain
\begin{align*}
&\E[(V(\bar X_{T\w\tau},T\w\tau)-V(\bar Y_{T\w\tau},T\w\tau))^2]
- (V(\bar x,0)-V(\bar y,0))^2
\\
&\ge
-c\,\E\int_0^{T\w\tau}\Big[(V(\bar X_t,t)-V(\bar Y_t,t))^2+|\bar X_t-\bar Y_t|^2 \Big]dt
\\
&\ge
-c\,\E\int_0^T\Big[(V(\bar X_{t\w\tau},t\w\tau)-V(\bar Y_{t\w\tau},t\w\tau))^2+|\bar X_{t\w\tau}-\bar Y_{t\w\tau}|^2 \Big]dt.
\end{align*}
Denoting
\[
\mu_t=\E[(V(\bar X_{t\w\tau},t\w\tau)-V(\bar Y_{t\w\tau},t\w\tau))^2],
\qquad
\gamma_t=\E\int_0^t|\bar X_{s\w\tau}-\bar Y_{s\w\tau}|^2ds,
\qquad 0\le t\le T,
\]
the above reads
\[
\mu_T-\mu_0\ge-c\int_0^T\mu_tdt-c\gamma_T.
\]
Letting $\tilde\mu_t=\mu_{T-t}$,
\[
\tilde\mu_T\le\tilde\mu_0+c\int_0^T\tilde\mu_tdt+c\gamma_T.
\]
Since this inequality holds for every $0<T\le T_0$, Gronwall's lemma gives,
for $0<T\le T_0$,
\[
(U(\bar x,T)-U(\bar y,T))^2=\tilde\mu_T\le(\tilde\mu_0+c\gamma_T)e^{cT}
=[\E(U_0(\bar X_{T\w\tau})-U_0(\bar Y_{T\w\tau}))^2+c\gamma_T]e^{cT}.
\]
Recalling that $c$ depends only on $T_0$, the result follows.
\qed

\begin{lemma}\label{lem02}
(a)
Let $(\bar X,\bar Y)$ be a solution to the system
\[
d\bar X_t=-b(\bar X_t)dt+\sqrt{2}dW_t,
\qquad
d\bar Y_t=-b(\bar Y_t)dt+\sqrt{2}dW_t,
\qquad
(\bar X_0,\bar Y_0)=(x,y),
\]
driven by the same $d$-dimensional standard Brownian motion $W$. Then $|\bar X_t-\bar Y_t|
\le e^{Ct}|x-y|$ for all $t$, where $C$ is the Lipschitz constant
of $b(\cdot)$.

\noi (b)
For each $\eps>0$, let $(\bar X^{(\eps)},\bar Y^{(\eps)})$ be
the coupling from \eqref{c3} starting at $(x^{(\eps)},y^{(\eps)})$. Assume $(x^{(\eps)},y^{(\eps)})
\to(x,y)$ as $\eps\to0$. Then $(\bar X^{(\eps)},\bar Y^{(\eps)})\To(\bar X,\bar Y)$
in $D(\R_+,\R^{2d})$. Consequently, for every $t,\del>0$,
\[
\lim_\eps\PP\Big(\sup_{s\in[0,t]}|\bar X^{(\eps)}_s-\bar Y^{(\eps)}_s|>|x-y|e^{Ct}+\del\Big)=0.
\]
\end{lemma}

\proof
(a) Immediate from
$\bar X_t-\bar Y_t=x-y-\int_0^t(b(\bar X_s)-b(\bar Y_s))ds$ and the Lipschitz
continuity of $b$.

(b)
The generator of $\bar X$ is
\[
\bar\calL_\eps f(x)=\sum_{i\le2d}\rho_{\eps,i}(x)\pl_if(x)=\sum_{i\le2d}r_{\eps,i^*}(x+\eps k_i)\pl_if(x)
=\sum_{i\le2d}(\eps^{-2}+\eps^{-1}q_{\eps,i^*}(x+\eps k_i))\pl_if(x).
\]
Therefore the asserted weak convergence will follow once the hypotheses of \cite[Theorem IX.4.21 p.\ 558]{jacshi} are verified. In particular, defining
\[
\tilde b_\eps(x)=\sum_{i\le2d} q_{\eps,i^*}(x+\eps k_i)k_i
=\sum_{i\le2d}q_{\eps,i}(x-\eps k_i)(-k_i),
\]
one only needs to verify that $\tilde b_\eps\to-b$ locally uniformly. But
\[
\tilde b_\eps(x)+b_\eps(x)
=\sum_{i\le2d}[-q_{\eps,i}(x-\eps k_i)+q_{\eps,i}(x)]k_i,
\]
hence by \eqref{c01}, $|\tilde b_\eps+b_\eps|_\iy\le2dC\eps$. Thus the uniform convergence
$b_\eps\to b$ stated in \eqref{c0} implies
that of $\tilde b\to-b$, and weak convergence follows.

The last assertion of (b) is immediate from (a) and the weak convergence.
\qed

\noi{\bf Proof of Proposition \ref{prop1}.}
Fix $T$ and denote
\[
\kappa(\eta)=\limsup_{\eps\to0}\sup_{t\in(0,T]}\om(U^{(\eps)}(\cdot,t),\eta),
\qquad \eta>0.
\]
We will show that $\kap(0+)=0$ by bounding it from above. To this end, given $\eta$, assume w.l.o.g.\ that $\kap(\eta)>0$. Fix sequences $t_n\in(0,T]$, $\eps_n\to0$ and $(x_n,y_n)\in\calS^{(\eps_n)}$ such that $|x_n-y_n|\le\eta$
and
\[
|U^{(\eps_n)}(x_n,t_n)-U^{(\eps_n)}(y_n,t_n)|\to\kap(\eta).
\]

Let $(\bar X_n,\bar Y_n)$ denote the coupling \eqref{c3} starting at $(x_n,y_n)$ and with $\eps=\eps_n$. In Lemma \ref{lem01}, let $\tau=\tau_n=\inf\{t:|\bar X_n(t)-\bar Y_n(t)|\ge1\}$.
Denoting by $c$ an upper bound on $|U_{0,n}|_\iy$ for all $n$,
we have for all large $n$,
\begin{align*}
\frac{\kappa(\eta)^2}{2}&<(U^{(\eps_n)}(x_n,t_n)-U^{(\eps_n)}(y_n,t_n))^2
\le c_T\E[(U^{(\eps_n)}_0(\bar X_n(t_n\w\tau_n))-U^{(\eps_n)}_0(\bar Y_n(t_n\w\tau_n))^2] + c_T\gamma_n,
\end{align*}
where
\[
\gamma_n=\E\int_0^T|\bar X_n(s\w\tau_n)-\bar Y_n(s\w\tau_n)|^2ds.
\]
By Assumption \ref{assn2}, denoting by $\om_0$ the modulus of continuity of $u_0$,
\[
\om(U^{(\eps)}_0,\eta)\le\om_0(\eta)+\del(\eps),
\]
where $\del(0+)=0$.
Hence
\[
\kappa(\eta)^2\le 2c_T\E\{[\om_0(|\bar X_n(t_n\w\tau_n)-\bar Y_n(t_n\w\tau_n)|)+\del(\eps_n)]^2\}
+2c_T\gamma_n.
\]
By Lemma \ref{lem02}, with $c_1=2e^{CT}$,
\[
\lim_n\PP(\sup_{s\in[0,T]}|\bar X_n(s)-\bar Y_n(s)|>c_1\eta)=0.
\]
Hence
\[
\limsup_n\E\{[\om_0(|\bar X_n(t_n\w\tau_n)-\bar Y_n(t_n\w\tau_n)|)+\del(\eps_n)]^2\} \le \om_0(c_1\eta)^2.
\]

As for a bound on $\gamma_n$, note that the processes $\bar X_n$ and $\bar Y_n$ have jumps of size $\eps_n$. Hence $|\bar X_n(\cdot\w\tau_n)-\bar Y_n(\cdot\w\tau_n)|\le 1+2\eps_n<\sqrt2$ for all large $n$. Thus the random variables
\[
Z_n=\int_0^T|\bar X_n(s\w\tau_n)-\bar Y_n(s\w\tau_n)|^2ds
\]
are bounded by $2T$ for all large $n$, whereas $(Z_n-T(c_1\eta)^2)^+\to0$ in probability. This shows that
\[
\limsup_n\gamma_n\le T(c_1\eta)^2.
\]
Combining these estimates, we get $\kappa(\eta)^2\le 2c_T\om_0(c_1\eta)^2+2c_TT(c_1\eta)^2$. This bound implies that $\kappa(0+)=0$, and the result follows.
\qed

\section{From ODE to PDE}\label{sec4}

Consider a solution $(\bu^{(\eps)},\bLa^{(\eps)})$ of \eqref{01} and construct from it the (deterministic) objects
\[
\xi^{(\eps)}_t(dx)=\sum_{y\in\calS_\eps}\bu^{(\eps)}(y,t)\del_y(dx),
\]
\[
\beta^{(\eps)}(dx,dt)=\sum_{y\in\calS_\eps}\bLa^{(\eps)}(y,t)\del_y(dx)dt.
\]
In this section we prove the following.
\begin{proposition}\label{prop4}
The family $(\xi^{(\eps)},\beta^{(\eps)})_{\eps>0}$ is relatively compact in $C(\R_+,\calM_1(\R^d))\times\calM^{(1)}(\R^d\times\R_+)$. Moreover, for any limit point $(\xi,\beta)$ as $\eps\to0$, $\xi$ is given by $\xi_t(dx)=u(x,t)dx$, $u\in\U(\R^d\times\R_+)$, and $(u,\beta)$ is a solution to \eqref{c1} (in the sense of Definition \ref{def1}).
\end{proposition}

\proof
Step 1: Relative compactness.

Let $\bar\bu^{(\eps)}$ be the solution to a version of \eqref{01} without the removal term, namely
\[
\bar\bu(x,t)=\bu_0(x)+\int_0^t(\calL^*_\eps \bar\bu+\bar\bu)(x,s)ds,
\quad(x,t)\in\calS_\eps\times\R_+.
\]
Then $\bu^{(\eps)}(x,t)\le\bar\bu^{(\eps)}(x,t)$ for all $x,t$. Moreover, $\bar\bu^{(\eps)}(x,t)=e^t\hat\bu^{(\eps)}(x,t)$, where $\hat\bu^{(\eps)}$ is the solution to
\begin{equation}\label{q1}
\hat\bu(x,t)=\bu_0(x)+\int_0^t\calL^*_\eps\hat\bu(x,s)ds,
\quad(x,t)\in\calS_\eps\times\R_+.
\end{equation}
Now let $X^{(\eps)}$ be a Markov process on $\calS_\eps$ with initial distribution $\bu_0(x)dx$ and generator $\calL_\eps$. Then \eqref{q1} is the forward Kolmogorov equation for this process, and therefore
\[
\bar\bu^{(\eps)}(\calB_r^c,t)\le e^t\hat\bu^{(\eps)}(\calB_r^c,t)
=e^t\,\PP(X^{(\eps)}_t\in\calB_r^c).
\]
A proof along the lines of that of Lemma \ref{lem02}(b) shows that $X^{(\eps)}\To X$ in $D(\R_+,\R^d)$ as $\eps\to0$, where $X$ is a solution to
\[
dX_t=b(X_t)dt+\sqrt{2}dW_t, \qquad X_0\sim u_0(x)dx.
\]
It follows that for every $T>0$, $\del>0$ there exists $r>0$ such that
\begin{equation}\label{q2}
\limsup_{\eps\to0}\sup_{t\in[0,T]}\bu^{(\eps)}(\calB_r^c,t)\le\del.
\end{equation}

Next, let $\ph\in C^\iy_c(\R^d)$. With $\lan\cdot,\cdot\ran=\lan\cdot,\cdot\ran_{\calS_\eps}$, we can write for the solution $(\bu^{(\eps)},\bLa^{(\eps)})$ of \eqref{01}, for a.e.\ $t$,
\begin{align}\label{q5}
\notag
\frac{d}{dt}\lan\ph,\bu^{(\eps)}\ran(t)
&=\lan\ph,\calL^*_\eps\bu^{(\eps)}+\bu^{(\eps)}-\bLa^{(\eps)}\ran(t)
\\
&=\lan\calL_\eps\ph+\ph,\bu^{(\eps)}\ran(t)-\lan\ph,\bLa^{(\eps)}\ran(t)
\end{align}
where \eqref{q4} is used.
For $x\in\calS_\eps$,
\begin{align*}
\calL_\eps\ph(x) &= \sum_{i\le 2d}r_{\eps,i}(x)(\ph(x+\eps k_i)-\ph(x))
\\
&=\sum_{i\le 2d}(\eps^{-2}+\eps^{-1}q_{\eps,i}(x))(\ph(x+\eps k_i)-\ph(x))=A_\eps(x)+B_\eps(x),
\end{align*}
where, denoting partial derivatives by $\ph_i$, $\ph_{ij}$,
\begin{align*}
A_\eps(x)&:= \eps^{-2}\sum_{i\le d}(\ph(x+\eps k_i)+\ph(x-\eps k_i)-2\ph(x))
\\
&=\eps^{-2}\sum_{i\le d}\Big(\ph_i(x)\eps+\frac{1}{2}\ph_{ii}(y_i)\eps^2
-\ph_i(x)\eps+\frac{1}{2}\ph_{ii}(\hat y_i)\eps^2\Big),
\\
B_\eps(x)&:= \eps^{-1}\sum_{i\le d}(q_{\eps,i}(x)(\ph(x+\eps k_i)-\ph(x))+q_{\eps,i^*}(x)(\ph(x-\eps k_i)-\ph(x)) )
\\
&=
\eps^{-1}\sum_{i\le d} (q_{\eps,i}(x)\ph_i(z_i)\eps-q_{\eps,i^*}(x)\ph_i(\hat z_i)\eps).
\end{align*}
Above, $y_i$ and $z_i$ (respectively, $\hat y_i$ and $\hat z_i$) are points on the segment connecting $x$ and $x+\eps k_i$ (respectively, $x-\eps k_i$). Because $\ph$ is smooth and compactly supported, there exists $c$ depending only on $\ph$ such that $|\ph_i(y_i)-\ph_i(x)|+|\ph_{ii}(y_i)-\ph_{ii}(x)|\le c\eps$. A similar statement holds for the remaining intermediate points. Hence
\[
A_\eps(x)=\Del\ph(x)+O(\eps)
\]
\[
B_\eps(x)=\sum_{i\le d}(q_{\eps,i}(x)-q_{\eps,i^*}(x))\ph_i(x)+O(\eps)
=b_\eps(x)\cdot\grad\ph(x)+O(\eps),
\]
where the constants in the expressions $O(\eps)$ depend only on $\ph$, and we have used \eqref{22b}. In view of Assumption \ref{assn1}, this shows
\[
\lim_{\eps\to0}\sup_{x\in\calS_\eps}|\calL_\eps\ph(x)-\calL\ph(x)|=0.
\]
Going back to \eqref{q5}, and using $|\bu^{(\eps)}(\cdot,t)|_1=1$, we have
\begin{equation}\label{q6}
\lan\ph,\bu^{(\eps)}\ran(t)-\lan\ph,\bu^{(\eps)}_0\ran
=\int_0^t\lan\calL\ph+\ph,\bu^{(\eps)}\ran(s)ds-\int_0^t\lan\ph,\bLa^{(\eps)}\ran(s)ds+\nu_\eps(t),
\end{equation}
where $\nu_\eps(t)\to0$ as $\eps\to0$ uniformly in $t$.

We use the above to show relative compactness of the measures $\beta^{(\eps)}$.
Fix $T$. Given $\del>0$, use \eqref{q2} to find $r$ such that, for all sufficiently small $\eps$, $\sup_{t\in[0,T]}\bu^{(\eps)}(\calB_r^c,t)<\del$. Let $\ph\in C^\iy_c(\R^d)$ be equal to $1$ in $\calB_r$ and supported in $\calB_{2r}$. Then by \eqref{q6}, and using the fact that $\calL\ph+\ph=1$ in $\calB_r$,
\begin{align*}
\beta^{(\eps)}(\calB_r\times[0,T])
&\ge -|\ph|_\iy\bu^{(\eps)}(\calB_r^c,T)
+\int_0^T\bu^{(\eps)}(\calB_r,s)ds
- |\calL\ph+\ph|_\iy \int_0^T \bu^{(\eps)}(\calB_r^c,s)ds+\nu_\eps(T)
\\
&\ge -c\del+T+\nu_\eps(T),
\end{align*}
for $c$ that depends only on $\ph$ and $T$. Since $\beta^{(\eps)}(\R^d\times[0,T])=T$, we have shown that for every $T>0$, $\del>0$ there exists $r$ such that
\[
\limsup_{\eps\to0}\beta^{(\eps)}(\calB_r^c\times[0,T]) < \del.
\]
This shows that the family of measures $\beta^{(\eps)}$, $\eps>0$ is relatively compact in $\calM^{(1)}(\R^d\times\R_+)$.

Next, to show that $\xi^{(\eps)}$ are relatively compact in $C(\R_+,\calM_1(\R^d))$, it is required to show that for every $T>0$ and $\del>0$ there is $r$ such that
\begin{equation}\label{z1}
\limsup_{\eps\to0}\sup_{t\in[0,T]}\xi_t^{(\eps)}(\calB_r^c)<\del,
\end{equation}
and for every $T>0$ and $\del>0$ there is $\eta>0$ such that
\begin{equation}\label{z2}
\limsup_{\eps\to0}w_T(\xi^{(\eps)}_\cdot,\eta)<\del.
\end{equation}
Here, the modulus of continuity $w$ is defined w.r.t.\ a metric on $\calM_1(\R^d)$ that is compatible with weak convergence. We take it to be the Levy--Prohorov metric
\[
\rho(\mu,\nu)=\inf\{\del>0:\mu(F)\le \nu(F^\del)+\del \text{ for all closed sets $F\subset\R^d$}\},
\]
where $F^\del$ is the $\del$-neighborhood of $F$ in $(\R^d,|\cdot|)$.

Note that \eqref{z1} is precisely \eqref{q2}. It remains to prove \eqref{z2}. To this end, let $F\subset\R^d$ be a closed set and let $s,t\in[0,T]$ be such that $0<t-s\le\eta$ for some $\eta>0$. By \eqref{t0}, for any $\del>0$,
\begin{align*}
\bu^{(\eps)}(F^\del,t)&=\sum_{x\in\calS_\eps}\bu^{(\eps)}(x,s)\frs_{t-s}(x,F^\del)-\int_s^t\sum_{x\in\calS_\eps}\bLa^{(\eps)}(x,\tau)\frs_{t-\tau}(x,F^\del)d\tau
\\
&\ge
\sum_{x\in\calS_\eps\cap F}\bu^{(\eps)}(x,s)\frp^{(\eps)}_{t-s}(x,F^\del)-e^{t-s}(t-s)
\\
&\ge
\bu^{(\eps)}(F,s)-\sum_{x\in\calS_\eps\cap F}\bu^{(\eps)}(x,s)\sum_{y\in\calB_\del^c(x)}\frp^{(\eps)}_{t-s}(x,y)-2\eta
\\
&\ge
\bu^{(\eps)}(F,s)-\sup_x\frp^{(\eps)}_{t-s}(x,\calB_\del^c)-2\eta.
\end{align*}
Now, by uniform ellipticity and boundedness of $b$, given $\del>0$, $\del'>0$, one can find $\eta>0$ such that $\frp_t(x,\calB_\del^c)<\del'$ for all $(x,t)\in\R^d\times(0,\eta]$. Therefore the same holds for $\frp^{(\eps)}$ provided $\eps$ is small. Thus, given $\del,\del'>0$ we can find $\eta$ such that the second term is bounded by $\del'$. If we now choose $\del'=\del/2$ and $\eta<\del/4$, we get
\[
\bu^{(\eps)}(F^\del,t)\ge\bu^{(\eps)}(F,s)-\del.
\]
Since this holds for all $F$, we get $\rho(\xi_s,\xi_t)\le\del$ provided $\eta$ is small. This proves \eqref{z2}.

Step 2:
Let $(\xi^{(\eps_n)},\beta^{(\eps_n)})$, $\eps_n\to0$ be a convergent sequence and denote its limit by $(\xi,\beta)$. We will prove now that
\begin{equation}\label{r1}
\text{there exists $u\in\U(\R^d\times\R_+)$ such that $\xi_t(dx)=u(x,t)dx$.}
\end{equation}
Fix $T>0$. It follows from the uniform estimate on the modulus of continuity stated in Proposition \ref{prop1} that there exists, for every $n$, an extension of $U^{(\eps_n)}$ from $\calS_{\eps_n}\times[0,T]$ to all of $\R^d\times[0,T]$, and there exists $\hat\om:\R_+\to\R_+$, $\hat\om(0+)=0$, such that
\[
|U^{(\eps_n)}(x,t)-U^{(\eps_n)}(y,t)|\le\hat\om(|x-y|),
\qquad x,y\in\R^d,\,t\in[0,T],\,n\in\N.
\]
In view of Lemma \ref{lem01}(a), we also have that $U^{(\eps_n)}$ are uniformly bounded on $\R^d\times[0,T]$. Let $\Q_T=[0,T]\cap\Q$, where $\Q$ denotes the set of rational numbers. Thus we can extract a subsequence $n_k\to\iy$ as $k\to\iy$, and a function $u:\R^d\times\Q_T\to\R_+$ such that $U^{(\eps_{n_k})}\to u$ uniformly in $\R^d\times\Q_T$, as $k\to\iy$, where the limit $u$ has $\hat\om$ as its modulus of continuity in $x$ for every $t\in\Q_T$.
(Note that, because we have not proved continuity of $t\mapsto U^{(\eps)}(x,t)$ uniformly in $\eps$, we cannot at this stage select $u$ to be continuous in $t\in[0,T]$).
Let $\ph\in C^\iy_c(\R^d)$ be supported in $\calB_r$, some $r>0$. For $y\in\calS_\eps$, denote the $\eps$-cube by $Q_\eps(y)=\prod_{i\le d}[y_i,y_i+\eps)$. Abbreviating $\eps_{n_k}$ to $\eps$, for fixed $t\in\Q_T$,
\begin{align*}
\lan\ph,\xi^{(\eps)}_t\ran
&=\sum_{y\in\calS_\eps}\ph(y)U^{(\eps)}(y,t)\eps^d
=\sum_{y\in\calS_\eps}\ph(y)U^{(\eps)}(y,t)\int_{Q_\eps(y)}dz
\\
&=
\sum_{y\in\calS_\eps}\int_{Q_\eps(y)}
\{\ph(z)U^{(\eps)}(z,t)+[\ph(y)U^{(\eps)}(y,t)-\ph(z)U^{(\eps)}(z,t)]\}dz.
\end{align*}
The expression in square brackets above is bounded by $|\ph|_\iy\hat\om(d^{1/2}\eps)+|U^{(\eps)}|_\iy\om_\ph(d^{1/2}\eps)$ whenever $z\in Q_\eps(y)$.
Moreover, this expression vanishes for $z$ outside $\calB_{2r}$ provided $\eps$ is small. Hence
\[
\lan\ph,\xi^{(\eps_{n_k})}_t\ran=\int_{\R^d}\ph(z)U^{(\eps_{n_k})}(z,t)dz+\bar\nu_k
\]
for some $\bar\nu_k\to0$ as $k\to\iy$. By the uniform convergence $U^{(\eps_{n_k})}\to u$, this shows that $\lan\ph,\xi^{(\eps_{n_k})}_t\ran\to\int\ph(z)u(z,t)dz$. However, we know that the full sequence $\xi^{(\eps_n)}_t$ converges to $\xi_t$ in $\calM(\R^d)$ as $n\to\iy$. This implies that $\xi_t(dx)=u(x,t)dx$. This conclusion holds for all $t\in\Q_T$. It remains to argue that $u$ has an extension to $C(\R^d\times[0,T],\R_+)$ such that $\xi_t=u(x,t)dx$ for all $t$.
This follows easily from continuity of $t\mapsto\xi_t$ in $\calM_1(\R^d)$ and the continuity in $x$ uniformly in $\R^d\times\Q_T$.

Step 3: $(u,\beta)$ solves \eqref{c1}.

Taking limits in \eqref{q6} and using $\xi_t=u(x,t)dx$ shows
\[
\lan\ph,u\ran(t)-\lan\ph,u_0\ran=\int_0^t\lan\calL\ph+\ph,u\ran(s)ds-\int_0^t\lan\ph,\beta_s\ran ds.
\]
It remains to show that $\beta\in\bB^u$.
By \eqref{01}, for all $n$,
\[
\sum_{x\in\calS_{\eps_n}}\int_0^T(|U^{(\eps_n)}(\cdot,t)|_\iy-U^{(\eps_n)}(x,t))\bLa^{(\eps_n)}(x,t)dt=0,
\]
thus
\[
\int_{\R^d\times[0,T]}(|U^{(\eps_n)}(\cdot,t)|_\iy-U^{(\eps_n)}(x,t))\beta^{(\eps_n)}(dx,dt)=0.
\]
Given $\del>0$ let $n_0$ be such that for all $n>n_0$, $|U^{(\eps_n)}-u|<\del T^{-1}$ on $\calS_{\eps_n}\times[0,T]$. Then for $n>n_0$,
\[
\int_{\R^d\times[0,T]}(|u(\cdot,t)|_\iy-u(x,t))\beta^{(\eps_n)}(dx,dt)<\del.
\]
Because $\beta^{(\eps_n)}\to\beta$, we get
\[
\int_{\R^d\times[0,T]}(|u(\cdot,t)|_\iy-u(x,t))\beta(dx,dt)\le\del.
\]
Sending $\del\to0$ we conclude that $\beta\in\bB^u$.
\qed

\section{PDE uniqueness}\label{sec5}

In this section we prove the following result.

\begin{proposition}\label{prop5}
There is at most one solution $(u,\beta)$ to \eqref{c1}  in $\U(\R^d\times\R_+)\times\calM^{(1)}(\R^d\times\R_+)$.
\end{proposition}

With this result, we have all that is needed in order to prove Theorem \ref{th1}.

\noi{\bf Proof of Theorem \ref{th1}.}
Part (a) follows from Proposition \ref{prop4} (existence) and Proposition~\ref{prop5} (uniqueness).

As for part (b), according to Lemma \ref{lem03} and Proposition \ref{prop2}, for every $\eps>0$ there exists a unique solution $(\bu^{(\eps)},\bLa^{(\eps)})$ to \eqref{01}. By Proposition \ref{prop2}, for fixed $\eps$, $(\xi^{\eps,N},\beta^{\eps,N})\to(\xi^{(\eps)},\beta^{(\eps)})$ in probability as $N\to\iy$, where
\[
\xi^{(\eps)}_t(dx)=\sum_{y\in\calS_\eps}\bu^{(\eps)}(x,t)\del_y(dx),
\qquad
\beta^{(\eps)}(dx,dt)=\sum_{y\in\calS_\eps}\bLa(y,t)\del_y(dx)dt.
\]
Moreover, combining Propositions \ref{prop4} and \ref{prop5} and denoting by $(\xi,\beta)$ the unique solution to \eqref{c1}, shows that $(\xi^{(\eps)},\beta^{(\eps)})\to(\xi,\beta)$ in $C(\R_+,\calM_1(\R^d))\times\calM^{(1)}(\R^d\times\R_+)$ as $\eps\to0$. The existence of $\eps_N\to0$ for which $(\xi^{\eps_N,N},\beta^{\eps_N,N})\to(\xi,\beta)$ in probability as $N\to\iy$ therefore follows by diagonalization.
\qed

Towards the proof of Proposition \ref{prop5} we shall need the following.
Let $\bar u_{u_0}$ and $\hat u_{u_0}$ denote the unique solution (in the class of smooth functions that are bounded in $\R^d\times[0,T]$ for every $T$) to
\[
\pl_t\bar u=\calL^*\bar u+\bar u,\qquad \bar u(\cdot,0)=u_0
\]
and, respectively,
\[
\pl_t\hat u=\calL^*\hat u,\qquad \hat u(\cdot,0)=u_0.
\]
Clearly, $\bar u_{u_0}(x,t)=e^t\hat u_{u_0}(x,t)$.

\begin{lemma}\label{lem21}
Every solution $(u,\beta)$ of \eqref{c1} satisfies $u(x,t)\le \bar u_{u_0}(x,t)$.
As a result, for every $T>0$, $\lim_{r\to\iy}\sup_{t\in[0,T]}\sup_{x\in\calB_r^c}u(x,t)=0$.
\end{lemma}

\proof
First we note that Definition \ref{def1} implies that for any $\ph\in C^\infty_c(\R^d, [0,T])$ and any $0\leq \tau< t\leq T$
\[
\lan \ph(\cdot,t),u(\cdot,t)\ran=\lan\ph(\cdot, 0),u(\cdot, 0)\ran
+\int_0^t\lan \calL\ph+\ph+\partial_s\ph,u\ran(s) ds
-\int_0^t\lan \ph,\beta_s\ran ds,
\qquad t>0.
\]
Moreover, since $u\in C_b(\R^d\times(0,\iy),\R_+)$, $\ph$ this equality extends to all $\ph\in C^\infty_0(\R^d\times [0,T])$ (i.e to functions decaying uniformly to $0$  at $\infty$ in $\R^d$). 

Next, note that $\bar u$ solves the equation
\[
\lan \ph(\cdot,t),\bar u(\cdot,t)\ran=\lan\ph(\cdot, 0),\bar u(\cdot, 0)\ran+\int_0^t\lan \calL\ph+\ph+ \partial_s\ph,\bar u\ran(s)ds,
\]
for any $t>0$ and test function $\ph\in C^\infty_0(\R^d\times [0,T])$. By Definition \ref{def1} and the above comment, we obtain for $w=u -\bar u$ (recalling $w(\cdot, 0)=0)$,
\[
\lan \ph(\cdot,t), w(\cdot,t)\ran=
\int_0^t\lan \calL\ph+\ph+ \partial_s\ph,w\ran(s) ds-\int_0^t\lan\ph,\beta_s\ran(s)ds, \ \ \beta\in\bB^u.
\]
Let $\ph\in C_0^\infty(\R^d\times [0, T])$ be a solution of the backward equation $\ph_s+ \calL \ph+\ph =0$ satisfying the terminal condition $\ph(\cdot, t)=\psi\in C_c^\infty (\R^d)$ in the interval $[0, t]$, and assume $\psi\geq 0$. Then
$\ph\geq 0$ in all of $\R^d\times[0,t]$ by the maximum principle. Thus
 \[
 \lan \psi, w(\cdot,t)\ran=
 -\int_0^t\lan\ph,\beta_s\ran(s) ds\leq 0.
 \]
This implies $w\leq 0$, which proves the first assertion of the lemma.

As for the tail estimate, if $\frp_t(x,y)$ denotes the fundamental solution to $\pl_t\hat u=\calL^*\hat u$, then $\hat u_{u_0}$ is given by $\hat u_{u_0}(y)=\int u_0(x)\frp_t(x,y)dx$. Now, as $u_0$ is uniformly continuous, it converges uniformly to zero as $|x|\to\iy$. Moreover, $\frp$ satisfies
Gaussian estimates uniformly in $\R^d\times[0,T]$. This and the bound $u\le\bar u_{u_0}=e^t\hat u_{u_0}$ easily imply the claimed tail estimate.
\qed

\noi{\bf Proof of Proposition \ref{prop5}.}
Let $(u,\beta)$ and $(v,\gamma)$ be two solutions in the indicated class.
Note that if $u=v$ then, by Definition \ref{def1}, for all $0\le s<t$, $\int_s^t\lan\ph,\beta_\theta\ran d\theta=\int_s^t\lan\ph,\gam_\theta\ran d\theta$ for all test functions $\ph$. This implies that for every rectangle $R=[a_1,b_1]\times\cdots\times[a_d,b_d]\times[s,t]$, one has $\beta(R)=\gamma(R)$. Hence $\beta=\gamma$. It thus suffices to prove that $u=v$.

Let $\eta_\del(x)=\del^{-d}\eta(x/\del)$
where $\eta$ is the standard mollifier in $\R^d$
(symmetric, nonnegative, supported in $\calB_1$, $\int\eta=1$). For $\psi:\R^d\to\R$ denote $\psi^\del=\eta_\del *\psi$. For $\psi:\R^d\to\R^d$, denote by $\psi^\del$ the componentwise convolution $(\eta_\del *\psi_i)_{i\in[d]}:\R^d\to\R^d$.
Let $\ph\in C^\iy_c(\R^d)$, and use $\ph^\del$ as a test function in Definition \ref{def1}. Note that (suppressing $t$),
\[
\lan\ph^\del, u\ran=\lan\ph,u^\del\ran,
\qquad
\lan\Del\ph^\del, u\ran=\lan\Del\ph,u^\del\ran=\lan\ph,\Del u^\del\ran.
\]
Let $\beta_t^\del(x,t)=(\eta_\del*\beta_t)(x)$ (a function).
By the antisymmetry of $\pl_i\eta_\del$ for every $i$,
\begin{align*}
\int b_i\pl_i(\ph*\eta)u\, dx&=\int b_iu (\ph*\pl_i\eta)\,dx
\\
&=\int\int b_i(x)u(x)\ph(y)\pl_i\eta(x-y)dydx
\\
&=-\int\int b_i(x)u(x)\ph(y)\pl_i\eta(y-x)dxdy
\\
&=-\int\ph ((b_iu) *\pl_i\eta)\,dx
\\
&=-\int\ph \pl_i(b_iu)^\del\,dx.
\end{align*}
Hence
\[
\lan b\cdot\nabla\ph^\del,u\ran = - \lan \ph,\nabla\cdot(bu)^\del\ran.
\]
Using these facts in Definition \ref{def1},
\[
\lan\ph,u^\del\ran(t)=\lan\ph,u_0^\del\ran+\int_0^t\lan\ph,\Del u^\del-\nabla\cdot(bu)^\del+u^\del\ran(s)ds-\int_0^t\lan\ph,\beta_s^\del\ran ds.
\]
Since $u^\del$ and its derivatives are continuous, it follows that
\[
u^\del(\cdot,t)=u^\del_0+\int_0^t(\Del u^\del-\nabla\cdot(bu)^\del+u^\del)(\cdot,s)ds-\int_0^t\beta^\del_sds.
\]
Let $w=u-v$ and $\theta_\del=(bw^\del)-(bw)^\del$. Then
\[
w^\del(\cdot,t)=\int_0^t(\calL^*w^\del+w^\del+\nabla\cdot\theta_\del)(\cdot,s)ds+\int_0^t(\gamma^\del_s-\beta^\del_s)ds.
\]
Thus
\begin{equation}\label{23}
\frac{1}{2}|w^\del|^2_2(t)
=\int_0^t\lan w^\del,\calL^*w^\del+w^\del+\nabla\cdot\theta_\del+\gamma^\del_s-\beta^\del_s\ran ds.
\end{equation}

 We now fix $T$ and proceed to bound the various terms on the r.h.s.\ of \eqref{23} for $t\in[0,T]$.
First, note that $\nabla w^\del$ satisfies the tail condition $\sup_{x\in\calB_r^c}|\nabla w^\del|(x,t)\to0$ as $r\to\iy$, uniformly in $t\in[0,T]$. Indeed, $\nabla w^\delta=\nabla \eta^\delta * w$ and  $w$  satisfies such a tail condition due to Lemma \ref{lem21}. As a result, integration by parts gives
$\lan w^\del,\Del w^\del\ran=-|\nabla w^\del|^2_2$.
Also,
$\lan w^\del,\nabla\cdot(bw^\del)\ran=-\int\nabla w^\del\cdot b w^\del=-\frac{1}{2}\int \nabla(w^\del)^2\cdot b=\frac{1}{2}\int(w^\del)^2\nabla\cdot b$.
This gives
\[
|\lan w^\del,\nabla\cdot(bw^\del)\ran|\le c|w^\del|^2_2.
\]
Also,
\[
\lan w^\del,\nabla\cdot\theta_\del\ran=-\int\nabla w^\del\cdot\theta_\del.
\]
For any $a>0$,
\[
\Big|\int \nabla w^\del\cdot \theta_\del \Big| \le |\nabla w^\del|_2 |\theta_\del|_2 \le a|\nabla w^\del|^2_2+a^{-1}|\theta_\del|^2_2.
\]
Combining the above estimates,
\begin{align*}
\frac{1}{2}|w^\del|^2_2(t) &\le -(1-a)\int_0^t|\nabla w^\del|^2_2(s)ds
+ c\int_0^t|w^\del|^2_2(s)ds
+\int_0^ta^{-1}|\theta_\del|^2_2(s)ds
\\
&\qquad+\int_0^t\lan w^\del,\gamma^\del_s-\beta^\del_s\ran ds.
\end{align*}
Because $b$ and $w$ are uniformly continuous in $\R^d$ and $\R^d\times[0,T]$, resp., we have that $|\theta_\del|_\iy=|(bw)^\del-bw^\del|_\iy\to0$ as $\del\to0$.
Since $u,v\in\U(\R^d\times\R_+)$ and $b$ is bounded, it follows that $|\theta_\del|_2\to0$ as $\del\to0$.
We can therefore find $a=a_\del\to0$ as $\del\to0$ such that one has $\kap_\del:=a_\del^{-1}\int_0^T|\theta_\del|^2_2(s)ds\to0$ as $\del\to0$.
This gives, for $t\in[0,T]$,
\[
\frac{1}{2}|w^\del|^2_2(t) \le
c\int_0^t|w^\del|^2_2(s)ds+\kap_\del
+\int_0^t\lan w^\del,\gamma^\del_s-\beta^\del_s\ran ds.
\]
Now
$$\int_0^t\lan w^\del,\gamma^\del_s-\beta^\del_s\ran ds=\int_0^t\lan w^{\del\del},\gamma_s-\beta_s\ran ds,$$
where $w^{\del\del}$ is the second convolution
$w^{\del\del}:= \eta_\del*\eta_\del*w$.
Since, by the definition of $\U(\R^d\times\R_+)$, $x\mapsto w(x,t)$ is continuous in $x$ uniformly in $(x,t)\in\R^d\times[0,T]$, there exists $\alpha_\delta$ such that $\sup_{t\in[0,T]}|w-w^{\del\del}|_\infty(t) <\alpha_\del$ and $\lim_{\del\rightarrow 0}\alpha_\del=0$.  It follows that
$$\int_0^t\lan w^{\del\del},\gamma_s-\beta_s\ran ds \leq t\alpha_\del+ \int_0^t\lan w,\gamma_s-\beta_s\ran ds\leq t\alpha_\del \ . $$
The last inequality follows by the monotonicity of the set valued operator $u\rightarrow \bB^u$, namely $\lan \beta_s-\gamma_s, u(\cdot,s)-v(\cdot,s)\ran \geq 0$ for a.e.\ $s$.  Indeed, 
$\lan \beta_s,u(\cdot,s)\ran = |u(\cdot,s)|_\infty$, 
$\lan \beta_s,v(\cdot,s)\ran \leq |v(\cdot,s)|_\infty$
for any $\beta\in \bB^u$ and a.e.\ $s$, by definition (\ref{betau}). 
Thus for all small $\del$,
\[
\frac{1}{2}|w^\del|^2_2(t) \le
c\int_0^t|w^\del|^2_2(s)ds +\kap_\del+t\alpha_\delta,
\qquad t\in[0,T].
\]
Using Gronwall's lemma, sending $\del\to0$ ,
\[
\lim_{\del\to0}\sup_{t\in[0,T]}|w^\del(\cdot,t)|_2=0.
\]
By the continuity of $w$, this shows that $w=0$ on $\R^d\times[0,T]$, and the proof is complete.
\qed

\appendix

\section{About Assumption \ref{assn1}}\label{app1}

Here we show that given $b\in C^1_b(\R^d,\R^d)$ with $b$ and $\nabla b$ globally Lipschitz, there always exists $q_\eps$ such that $(q_\eps,b_\eps,b)$ satisfy Assumption \ref{assn1}.

For example, set $g_i=(b\cdot e_i)^+*\eta$, $1\le i\le d$,
where $\eta$ is the standard mollifier in $\R$,
and $g_{i+d}=-b\cdot e_i+g_i$.
Take $q_{\eps,i}(x)=g_i(x)$ for $1\le i\le 2d$ and $x\in\calS_\eps$.
Now, $g_i\ge0$ for $1\le i\le d$, whereas $g_i$ are bounded below for $d+1\le i\le2d$. As a result,
$\min_{i\le 2d}\inf_xr_{\eps,i}(x)>0$ provided $\eps$ is sufficiently small.
Finally, $g_i$ and their derivatives are Lipschitz
for all $1\le i\le2d$ by the assumed properties of $b$,
and thus \eqref{c0} and \eqref{c01} follow.

{\bf Acknowledgment.} RA was supported in part by ISF 3240/25. LM was supported in part by ISF 1985/22.

%\footnotesize

%%\bibliographystyle{plain}
%%\bibliographystyle{annotate}
%%\bibliographystyle{apalike}
\bibliographystyle{is-abbrv}

\bibliography{main}

\vspace{.5em}

\end{document}